\def\VERSION{29.1.2023}
\def\WHO{nbd} 
\def\users{us}    
\def\users{world} 
\documentclass[12pt]{article} 
\usepackage{color}
\usepackage{amsmath}
\usepackage{amsthm}
\usepackage{amssymb}
\usepackage{epsfig}
\usepackage{psfrag}
\usepackage{graphicx}
\usepackage{textcomp}
\usepackage{bm}
\usepackage{pgf}
\numberwithin{equation}{section}
\usepackage{upgreek} 
\newtheorem{theorem}{Theorem}[section]

\newtheorem{definition}[theorem]{Definition}

\newtheorem{proposition}[theorem]{Proposition}

\newtheorem{remark}[theorem]{Remark}

\usepackage{mathrsfs,cite}
\usepackage{ifthen}
\usepackage{ulem}
\usepackage[colorlinks=true, pdfstartview=FitV, linkcolor=blue, 
            citecolor=blue, urlcolor=blue]{hyperref}
\usepackage{cancel}\ifthenelse{\equal{\users}{world}}
{
\newcommand{\REM}[1]{}

	\newcommand{\DELETE}[1]{}

        \newcommand{\COMMENT}[1]{}
        \newcommand{\TCOMMENT}[1]{}

    \newcommand{\MARGINOTE}[1]{}
}	
{
\definecolor{brown}{rgb}{0.6,0.2,0.2}
\newcommand{\REM}[1]{\marginpar{\bfseries\tiny{\color{blue}#1}}}

 \newcommand{\COMMENT}[1]{{\color{blue}\uuline{#1}\color{black}}} 
 \newcommand{\DELETE}[1]{{\color{brown}\cancel{#1}\color{black}}}

 \newcommand{\TCOMMENT}[1]{{\color{blue}{ #1}}}

\newcommand{\MARGINOTE}[1]{\marginpar{\color{red}\tiny\texttt{#1}}}
\usepackage{fancyhdr}
\pagestyle{fancy}
\headheight=28pt\headwidth=17cm
\definecolor{gray}{gray}{0.5}
\newcount\hour \newcount\minute
\hour=\time
\divide \hour by 60
\minute=\time
\loop \ifnum \minute > 59 \advance \minute by -60 \repeat

\rhead{\color{gray}paleomagnetism\\
by T.Roub\'\i\v cek}
\chead{}
\lhead{Version \VERSION, file: \jobname.tex, \underline{\Large\color{red}\WHO's working on} \\
compiled:
\number\day.\number\month.\number\year\ at
\the\hour:\ifnum\minute<10 0\fi\the\minute\ h }
}

\newcommand{\R}{\mathbb{R}}

\newcommand{\bbI}{\mathbb{I}}

\newcommand{\bbD}{\mathbb{D}}
\newcommand{\bbC}{\mathbb{C}}
\newcommand\DT[1]{\mathchoice
                 {{\buildrel{\hspace*{.1em}\text{\LARGE.}}\over{#1}}}
                 {{\buildrel{\hspace*{.1em}\text{\LARGE.}}\over{#1}}}
                 {{\buildrel{\hspace*{.1em}\text{\Large.}}\over{#1}}}
                 {{\buildrel{\hspace*{.1em}\text{\large.}}\over{#1}}}}

\newcommand{\lineunder}[2]{\LU{\begin{array}[t]{c}\underbrace{#1}\vspace*{.5em}\end{array}}{\mbox{\footnotesize\rm #2}}}
\newcommand{\linesunder}[3]{\LSU{\begin{array}[t]{c}\underbrace{#1}\vspace*{.5em}\end{array}}{\mbox{\footnotesize\rm #2}}{\mbox{\footnotesize\rm#3}}}
\newcommand{\LU}[2]{\begin{array}[t]{c}#1\vspace*{-1em}\\_{#2}\end{array}}
\newcommand{\LSU}[3]{\begin{array}[t]{c}#1\vspace*{-1em}\\_{#2}\vspace*{-.5em}\\_{#3}\end{array}}
\renewcommand{\d}{{\rm d}}

\newcommand{\NablaS}{\Nabla_{\scriptscriptstyle\textrm{\hspace*{-.3em}S}}^{}}
\newcommand{\divS}{\mathrm{div}_{\scriptscriptstyle\textrm{\hspace*{-.1em}S}}^{}}

\def\vv{{\bm v}}

\def\nn{{\bm n}}
\def\ff{{\bm g}}
\def\mm{{\bm m}}
\def\bb{{\bm b}}
\def\hh{{\bm h}}
\newcommand{\RR}{\bm{R}}
\newcommand{\RRk}{\bm{R}_\etau^k}
\newcommand{\rr}{\bm{r}}
\newcommand{\rrk}{\bm{r}_\etau^k}
\newcommand{\GM}{M}
\newcommand{\KK}{K}
\newcommand{\DD}{\nu_1}
\newcommand{\W}{w}
\newcommand\EE{{\bm e}}
\newcommand{\rexp}{r}
\newcommand{\qexp}{2}

\def\etau{{\eps\tau}}

\def\vvk{\vv_\etau^k}
\def\vvkk{\vv_\etau^{k-1}}
\def\mmk{\mm_\etau^k}
\def\thetak{\theta_\etau^k}
\def\mmkk{\mm_\etau^{k-1}}
  \def\Eek{\Ee_\etau^k}
  \def\Eekk{\Ee_\etau^{k-1}}
 \def\Eetau{\Ee_\etau^{}}
\def\overlineEetau{\hspace*{.2em}\overline{\hspace*{-.2em}\bm E}_\etau^{}}
\def\overlineSetau{\hspace*{.2em}\overline{\hspace*{-.2em}\bm S}_{\text{\sc e},\etau^{}}}
\def\overlineSvtau{\hspace*{.2em}\overline{\hspace*{-.2em}\bm S}_{\text{\sc v},\etau^{}}}
\def\overlineSstrtau{\hspace*{.2em}\overline{\hspace*{-.2em}\bm S}_{\text{\sc c},\etau}^{}}

\def\overlineRRtau{\hspace*{.2em}\overline{\hspace*{-.2em}\RR}_{\etau}^{}}
\def\overlinemmtau{\hspace*{.15em}\overline{\hspace*{-.15em}\mm\hspace*{-.1em}}_{\etau}^{}}

\def\overlinevvtau{\hspace*{.15em}\overline{\hspace*{-.15em}\vv}_{\etau}^{}}

\def\NU{\nu_2}
\def\eps{\varepsilon}

\def\FF{\breve\varphi}
\def\GG{\widetilde\varphi}

\newcommand\ZJ[1]{\mathchoice
                 {{\buildrel{\hspace*{.1em}{_{\,\boldsymbol\circ}}}\over{#1}}}
                 {{\buildrel{\hspace*{.1em}{_{\,\boldsymbol\circ}}}\over{#1}}}
                 {{\buildrel{\hspace*{.1em}{\boldsymbol\circ}}\over{#1}}}
                 {{\buildrel{\hspace*{.1em}{\boldsymbol\circ}}\over{#1}}}}

\newcommand{\Nabla}{\nabla}

\def\Vdots{\!\mbox{\setlength{\unitlength}{1em}
\begin{picture}(0,0)
\put(-.07,0){.}
\put(-.07,.3){.}
\put(-.07,.6){.}
\end{picture}
}
}

\newcounter{myfigure}
\newenvironment{my-picture}[3]{\refstepcounter{myfigure}\label{#3}\setlength{\unitlength}{1em}\begin{picture}(#1,#2)}{\end{picture}}

\newcommand\DELETEDELETE[1]{}

\newcommand\pdt[1]{\frac{\partial{#1}}{\partial t}}
\newcommand\Ee{{\bm E}}              
\newcommand\Ep{{\bm\varPi}}              

\def\widetildeEp{\hspace*{.1em}\widetilde{\hspace*{-.1em}\Ep}}
\def\CC{\omega}
\def\CCC{\omega}

\begin{document}
\allowdisplaybreaks

\noindent{\LARGE\bf A thermodynamical model for paleomagnetism
\\[.2em] in Earth's crust. 
}

\bigskip\bigskip

\noindent{\large\sc Tom\'{a}\v{s} Roub\'\i\v{c}ek}\\
{\it Mathematical Institute, Charles University, \\Sokolovsk\'a 83,
CZ--186~75~Praha~8,  Czech Republic}\\and\\
{\it Institute of Thermomechanics, Czech Academy of Sciences,\\Dolej\v skova~5,
CZ--182~08 Praha 8, Czech Republic}

\bigskip\bigskip

\begin{center}\begin{minipage}[t]{16.5cm}

{\small

\noindent{\bfseries Abstract.}
\baselineskip=12pt
A thermodynamically consistent model for soft deformable viscoelastic magnets
is formulated in actual space (Eulerian) coordinates. The possibility of
a ferro-paramagnetic-type (or ferri-antiferromagnetic) transition exploiting
Landau phase transition theory as well as mechanical melting
or solidification is considered, being motivated and applicable to
paleomagnetism (involving both thermo- and isothermal and viscous
remanent magnetization) in rocks in Earth's crust and to rock-magma
transition. The temperature-dependent Jeffreys rheology in the deviatoric part
combined with the Kelvin-Voigt rheology in the spherical (volumetric) part
is used. The energy balance and the entropy imbalance behind the model are
demonstrated, and its analysis is performed by time discretization, proving
existence of weak solutions.

\medskip


\medskip


} 
\end{minipage}
\end{center}

\bigskip

\section{Introduction}

Magnetism in our planet Earth has two main mechanisms. The mains source
of geomagnetic field is in the fluidic part of the Earth core, called the
outer core, composed primarily from electrically conductive
hot Iron with Nickel flowing in high speed of the order 10\,km/year and
inducing magnetic field via magnetodynamo effect. The second mechanism
is in the very upper part of the silicate mantle, called the crust, which
is rather cold and exhibits para-to-ferro (or antiferro-to-ferrimagnetic
phase transition under the 
geomagnetic field generated primarily by magnetic dynamo in the fluidic
Iron-Nickel outer core, considered
time-space dependent but given in this paper. These two
very different magnetic phenomena are primarily related respectively with the
names of Hannes O.\,G.\,Alfv\'en and Louis E.\,F.\,N\'eel, both Nobel prize
winners in 1970. For completeness, let us mention that other contributions are
due to the (relatively weak) magnetodynamo effects in moving salty oceans or in
hot silicate mantle and due to the magnetic field from the Sun which interacts
mainly with the magnetosphere around the Earth and does not substantially
influence Earth's interior.

In this paper, we will focus on the second phenomenon of magnetism in
rocks containing magnetic materials, called {\it paleomagnetism}. Magnetism
in some rocks forms a vital part of rock physics and mechanics,
cf.\ \cite{Butl92PMDG,Camp03IGF,DunOzd97RM,LanMel06EM,StaBan74PPRM}. Magnetism
in the oceanic or the continental Earth crust is an important phenomenon
which gives information about history of geomagnetic field and about various
mechanical processes e.g.\ behind  folding and faulting of rocks or
formation mountains or even or even movement of continents in the far
history of the Earth or even in some other terestial-type planets/moons as
particularly Mars and ``our'' Moon.
As in all magnetic materials, the magnetism is heavily dependent on
temperature and, above certain ``critical'' temperature disappears because
spontaneous magnetism in ferromagnetic materials undergo the transition to a
non-magnetic paramagnetic variant. Ferromagnetism is caused by a
synchronized (parallel) orientation of spins of atoms in crystal lattices
while paramagnetism is characterized by random orientation. In the case
of ferro-to-para magnetic transition, the mentioned critical temperature
is called the {\it Curie temperature}.

Actually, in rocks, orientation of magnetic moments in neigbouring
atoms is opposite. Their magnitude can be the same or one of them 
can dominate. Then we speak about antiferro- or ferri-magnetism, respectively.
The spontaneous magnetization in antiferromagnets is zero, similarly
as in paramagnets. Like in ferro/para-magnetism, it depends on temperature
and, instead of the Curie temperature, we speak about the
{\it N\'eel temperature}.

There are several mechanisms behind natural remanent magnetism in cold
rocks. The most typical is related with the antiferro-to-ferri-magnetic
phase transition (phenomenologically quite similar to
para-to-ferro-magnetic phase transition) during cooling of initially hot
rock containing magnetic minerals in the geomagnetic field. This is
referred to as a {\it thermoremanent
magnetization}\index{magnetization!thermoremanent (TRM)}
(TRM), cf.~Figure~\ref{fig4}.
\begin{center}
\begin{my-picture}{38}{16}{fig4}
\hspace*{-0em}{\includegraphics[width=38em]{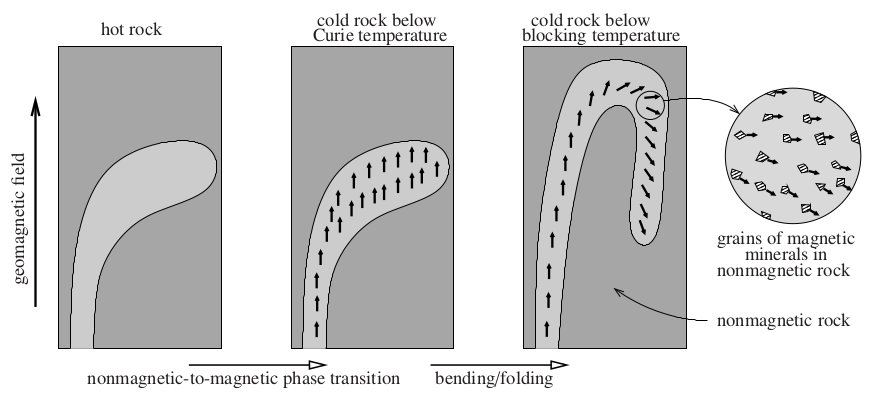}}
\end{my-picture}
\nopagebreak
\\
  {\small\sl\hspace*{-2em}Fig.~\ref{fig4}:~\begin{minipage}[t]{35em}
 A schematic illustration of thermoremanent magnetization in cooling and
 subsequent bending/folding of rocks in a constant geomagnetic field.
In fact, grains of magnetic minerals have randomly oriented easy-magnetization
axes, some of them being magnetized more while others less.
\end{minipage}
} 
\end{center}
Another magnetization in cold rock with magnetic minerals might be 
by strong magnetic fields typically due to lightning strikes, called an
{\it isothermal remanent magnetization} (IRM). Beside these two, there
is also a {\it viscous remanent magnetization} (VRM) which may
occur when rocks are exposed by a modern-day geomagnetic field which are
stronger than geomagnetic field but anyhow not so strong to lead to an
immediate (re)magnetization. Actually, there are some other processes
and mechanisms occasionally relevant in paleomagnetism, as chemical or
depositional (detrital) remanent magnetization but we will not
address them in our model. All these mechanisms belong to
{\it paleomagnetism}, cf.\ in particular \cite{CoxDoe60RP,Butl92PMDG}
and historical
references \cite{Irvi56PPAP,Neel55STAR}.

The scope of this paper is the following: In Section~\ref{sec-model}, we
briefly present relevant physics and formulate the governing system of partial
differential equations. Then, 
in Section~\ref{sec-thermodyn}, we justify the model as far as its desired
energy balances and the entropy imbalance. Last, in Section~\ref{sec-anal},
the analysis of an initial-boundary-value problem for the system from
Section~\ref{sec-model} is performed by time discretization, proving
existence of suitably defined weak solutions.

\section{The thermo-magneto-mechanical model}\label{sec-model}

Before writing the system of partial differential equations describing
the model, let us first articulate the main concepts and motivation
we want to employ, and related attributes of the model.

The main attribute of the model will be its formulation in {\it Eulerian
coordinates}. This reflects the reality that, within geological time
scales of instantaneous evolution of Earth, there is no reference
configuration. It contrasts to engineering work-pieces where a reference
configuration may refer to the shape in which they have been manufactured and
the Lagrangean description which uses such a reference configuartion is well
motivated even for largely strained elastic solids.

Rocks considered in a limited (and not much large) space-time region of
the crust can be considered not much compressible, and in particular,
the mass density does not vary substantially and can be (and often is)
considered constant. This simplifies the model and its analysis
substantially. Often, the models are simplified even more by imposing
incompressibility. Yet, incompressible models do not facilitate propagation
of longitudinal seismic waves, which substantially limit geophysical
applications, although in some situation (in particular in
paleomagnetism), even quasistatic models which neglects all dynamical
effects have reasonable applications. Yet, forgetting inertial
effects seems to bring analytical difficulties by lacking an
immediate control on acceleration. Therefore, we will use a
dynamical and not fully incompressible model but we will adopt a compromise
concept of only slightly compressible (so-called
{\it semi-compressible}) {\it materials} \cite{Roub19QSF} with mass
densities constant in space and time. Cf.\ Remark~\ref{compres} for a
fully compressible model.

Various gradient theories are applicable and allows for inventing
various dispersion of velocities of elastic waves into the model
and facilitate its analysis, cf.\ \cite{Roub19QSF}.
The particular enhancement by dissipative gradient terms exploits
the general ideas of {\it multipolar} (also called non-simple) media by
Green and Rivlin \cite{GreRiv64MCM} adopted for {\it fluids} by  
Ne\v cas at al.\ \cite{BeNeRa99EUFM,Neca94TMF,NeNoSi89GSIC,NecRuz92GSIV}
and Fried and Gurtin \cite{FriGur06TBBC}.
More specifically, we will use nonlinear 2nd-grade nonsimple fluids.

Essentially, rocks will then be considered as
viscoelastic fluids but very high viscosity (typically of the order
$10^{22\pm2}$ Pa\,s) unless being melted to magma when the viscosity
drops down substantially (typically being of the order
$10^{4\pm3}$ Pa\,s), cf.\ Figure~\ref{fig1}.
This rock-magma transition will be covered
by our model, too. Of course, the magma is hot and surely much above
the Curie temperature, thus nonmagnetic.

For thermoremanent magnetization of rocks,
we need a ferro-para-magnetic transition a formulated and
analyzed in \cite{PGRoTo10TCTF} by using the Landau \cite{Land37TPT} 
idea. Possibly, this may be quite equally
interpreted as ferro-antiferro-magnetic transition, too.

\begin{center}
\begin{my-picture}{35}{15}{fig1}
\hspace*{-0em}{\includegraphics[width=35em]{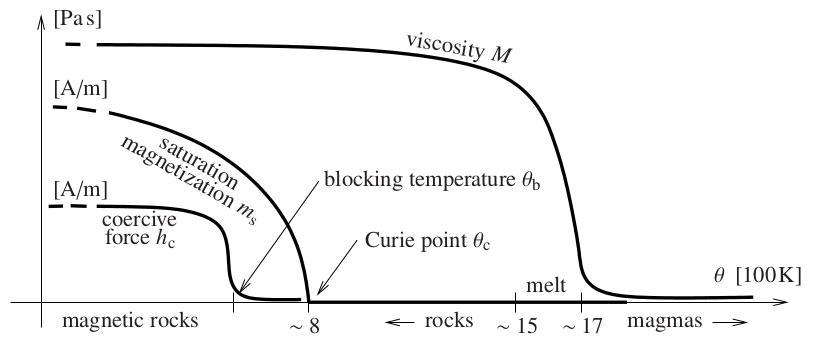}}
\end{my-picture}
\nopagebreak
\\
  {\small\sl\hspace*{-2em}Fig.~\ref{fig1}:~\begin{minipage}[t]{35em}
 Schematic temperature dependence of magnetic properties behind
 ferri-to-antiferro\-magnetic phase transition in solid rocks and
 of mechanical properties (viscosity) behind melting of rocks towards magma.
\end{minipage}
} 
\end{center}

Further important attribute is a proper choice of rates (time derivatives)
of intensive variables in the Eulerian setting. The adjective
``intensive'' refers to variables or properties which do not depend on
the system size or the amount of material in the system, as e.g.\
velocity or temperature, in contrast to extensive variables as e.g.\
momentum or entropy. Intensive scalar variables are transported by the
{\it convective} (also called {\it material}) {\it time derivative}
defined, for a scalar variable $\alpha$ and velocity field $\vv$, as
\begin{align}
\DT\alpha=\pdt\alpha+\vv{\cdot}\Nabla\alpha\,.
\end{align}
It is applicable componentwise also for a velocity vector itself for
which it reads as $\DT\vv=\pdt{}\vv+(\vv{\cdot}\Nabla)\vv$, i.e.\ the
last term means componentwise 
$[(\vv{\cdot}\Nabla)\vv]_i=\sum_{j=1}^d\vv_j\frac{\partial}{\partial x_j}\vv_i$.
This convective derivative is however not objective for general vector- or
tensor-valued intensive variables (as e.g.\ magnetization or stress)
which have to be transported by other time derivatives, however.
Objectivity here means that the time derivatives do not depend on the
frame of reference. For such rates, it is reasonable to require, beside
objectivity, also being so-called corotational.
The simplest corotational form is the {\it Zaremba-Jaumann time
derivative} \cite{Jaum11GSPC,Zare03FPTR} justified for stress
rates by M.\,Biot \cite[p.494]{Biot65MID}, cf.\ also \cite{Fial11GSSM},
defined as
\begin{align}\label{ZJ-tensor}
&\ZJ{\bm A}=
\pdt{\bm A}+(\vv{\cdot}\Nabla){\bm A}-{\bm W}{\bm A}+{\bm A}{\bm W}
\ \ \text{ with the material spin }\ 
{\bm W}={\rm skw}(\Nabla\vv)\,
\end{align}
with ``skw'' denoting the skew-symmetric part, i.e.\ ${\rm skw}(\Nabla\vv)
=\frac12(\Nabla\vv-(\Nabla\vv)^\top)$.
In isotropic materials, it is relevant also for strain rates, cf.\
\cite{Roub20SPTC}.
The important attribute of the corotational derivative is that it
commutes with the transposition and with convective derivatives for 
for traces and pressure-like tensors:
\begin{align}
&(\ZJ{\bm A})^\top\!\!=({\bm A}^\top)\!\ZJ{^{}}\,,\ \ \ \ \ \
{\rm tr}\ZJ{\bm A}=({\rm tr}{\bm A})\!\DT{^{}}\,,\ \ \text{ and }\ \ 
(\alpha\bbI)\!\ZJ{^{}}
=\DT\alpha\bbI\,.
\label{ZJ-property}\end{align}
These properties are important to keep properties of
symmetric or deviatoric preserved during evolution governed by
such derivatives as noticed e.g.\ in \cite{ThKaPo15LSRB}, and to gain the
expected convective derivative for transition from solid to fluid models
during melting of rocks to magmas as noticed in \cite{Roub20SPTC}.
Consistently with the Zaremba-Jaumann corotational derivative
used for stress or strain tensors in the model, it is natural to use it
also for vectors, cf.\ e.g.\ \cite[Sect.\,5.5]{Mart19PCM}, defined as
\begin{align}
\label{ZJ-vector}\ZJ\mm=\pdt\mm+(\vv{\cdot}\Nabla)\mm-{\bm W}\mm
\ \ \ \ \ \ \ \text{ where again}\ \ \ \ \ {\bm W}={\rm skw}(\Nabla\vv)\,.
\end{align}
Actually, it was used for magnetization in \cite[Sect.\,2.5]{BePGVa01DDFG}.

We neglect thermal expansion (and the related buoyancy effects --
cf.\ Remark~\ref{OB}) and also we neglect variation of mass density not only
by thermal expansion but also by pressure.
Anyhow, we admit the medium being slightly compressible. This is
a reasonable compromise for solid or fluidic geophysical materials
if one does not want to exclude pressure (longitudinal) waves like
it would happen in fully incompressible models. On the other hand,
allowing for a slight compressibility but taking fixed mass density
is relevant if pressure variations are negligible in comparison with
elastic bulk modulus (which is usually quite large of the order of GPa's).
It leads either to (slight) violation of energy
conservation or, when considering an extra (small) force invented by
Temam \cite{Tema69ASEN}, to (slight) violation of Galilean invariancy,
cf.\ Remark~\ref{compres}. Here, we will adopt the latter option. Fully
convective variants of such a compromise fluidic models are sometimes called
{\it quasi-compressible}, cf.\ \cite{Roub19QSF}. In the former
energy-violating variant, they are sometimes used in physical modelling,
cf.\ e.g.\ \cite{BATP20CFPF,PGYD20SHMM,SPGB05TMMD}.

We consider a fixed bounded domain $\varOmega\subset\R^d$ with $d=2$ or 3 and
with a Lipschitz boundary $\varGamma$ and a time interval $I=[0,T]$.
The basic variables and data for this preliminary model are summarized
in the following table:

{\small
\begin{center}
\begin{minipage}[t]{40em}
\hrule
\ \end{minipage}
\end{center}\vspace*{-2em}
\begin{center}
\ \ \ \ \begin{minipage}[t]{21em}

$\vv$ velocity (in m/s)

$\Ee$ the elastic strain tensor (symmetric)

$\Ep$ the inelastic strain tensor (deviatoric)

$\RR$ the inelastic strain rate (in s$^{-1}$)

$\mm$ magnetization vector (in A/m)

$\rr$ magnetization rate (in Am$^{-1}$s$^{-1}$)

$\bb$ magnetic induction (in T)

$\hh_{\rm drv}$ driving magnetic field (in A/m)

$\hh_{\rm dem}$ intensity of demagnetizing field (in A/m)

$u$ potential of demagnetizing field (in A)

$\theta$ temperature (in K)

$\W$ enthalpy (in Pa=J/m$^3$)

$\zeta$ magnetic dissipation potential (in Pa/s=W/m$^{3}$)\!\!\!\!\!\!

$\mu_0$ vacuum permeability ($\sim1.257{\times}10^{-6}$\,H/m)

\end{minipage}\hspace*{2em}
\begin{minipage}[t]{18em}

$\eta$ entropy (in Pa/K)

$\psi$ Helmholtz free energy (in Pa)

$\EE(\vv)\,{=}\,{\rm sym}(\Nabla\vv)$ small strain rate (in s$^{-1}$)

$\bm{S}_\text{\sc e}^{}$ elastic stress (symmetric - in Pa)

$\bm{S}_\text{\sc v}^{}$ viscous stress (symmetric - in Pa)

$\bm{S}_\text{\sc c}^{}$ capillarity/couple stress (in Pa)

$\mathfrak{H}$ hyperstress (in Pa\,m)

$\varrho$ mass density (in kg/m$^3$)
  
$\ff$ gravity acceleration (in ms$^{-2}$)

$\GM$ Maxwell viscosity modulus (in Pa\,s)

$\DD$ Stokes viscosity modulus (in Pa\,s)

$\KK$ heat conductivity (in Wm$^{-1}$K$^{-1}$) 

$j_{\rm ext}$ external heat flux (in W/m$^{2}$) 

$\hh_\text{\rm geo}$ intensity of geomagnetic field 

\end{minipage}
\end{center}
\vspace*{-1.5em}\begin{center}
\begin{minipage}[t]{40em}
\hrule
\ \end{minipage}
\end{center}
}

The basic ingredient of the model is a specific {\it free energy} $\psi$.
The simplest form  in terms of the elastic strain $\Ee$ with a rather weak
(linearized) coupling of mechanical and magnetic effects is
$\psi(\Ee,\mm,\theta)=\varphi_0(\Ee,\mm)+\theta\varphi_1(\mm)+
\mu_0\hh_{\rm dem}{\cdot}\mm-\phi(\theta)$ or, for the purpose of the analysis
below, rather as
\begin{align}\label{free}
\psi(\Ee,\mm,\theta)=\varphi(\Ee,\mm)
+\theta\CC(\mm)+\mu_0\hh_{\rm dem}{\cdot}\mm-\phi(\theta)\,,
\end{align}
where $\hh_{\rm dem}$ is the intensity of demagnetizing field
satisfying the equations
\begin{align}\label{Max}
{\rm curl}\,\hh_{\rm dem}=0\ \ \ \ \ \text{ and }\ \ \ \ \ {\rm div}(\hh_{\rm dem}
+\chi_{\varOmega}^{}\mm)=0
\ \ \ \ \text{ on }\ \ \R^d\,,
\end{align}
which is the quasistatic rest of the Maxwell electromagnetic system in
electrically non-conductive magnetic media. The second equation is the
Gauss law for magnetism ${\rm div}\,\bb=0$ with the magnetic induction
$\bb=\mu_0(\hh_{\rm dem}{+}\chi_\varOmega^{}\mm)$ with $\mu_0$ the vacuum
permeability. Here, $\chi_\varOmega^{}$ is the characteristic function
of the domain $\varOmega$, i.e.\ $\chi_\varOmega^{}(x)=1$ for $x\in\varOmega$
while =\,0 otherwise. The static system \eqref{Max} is thus
considered on the whole Universe $\R^d$ at each time instant.
The first equation in \eqref{Max} implies existence of
a (scalar) magnetic potential $u$ such that $\hh_{\rm dem}=-\Nabla u$.
Thus \eqref{Max} can be written as a single Poisson equation
\begin{align}\label{Max+}
\Delta u={\rm div}(\chi_\varOmega^{}\mm)\ \ \ \ \text{ on }\ \ \R^d
\end{align}
to be understood in the sense of distributions.

Let us now write the system of six partial differential equations (one of them
being an inclusion) for $\vv$, $\Ee$, $\Ep$, $\mm$, $u$, and $\theta$. More
specifically, it is composed from a momentum equation, the Green-Naghdi's
additive decomposition of the total strain written in terms of rates, a flow
rule for the inelastic strain $\Ep$, a flow rule for the magnetization $\mm$,
the Poisson equation for the demagnetizing-field potential $u$, and
the heat-transfer equation for temperature $\theta$. Namely,
\begin{subequations}\label{MG}\begin{align}\nonumber
&\varrho\DT\vv={\rm div}\big(\bm{S}_\text{\sc e}^{}
    +\bm{S}_\text{\sc v}^{}
+\bm{S}_\text{\sc c}^{}-\mu_0(\hh{\cdot}\mm)\bbI
-{\rm div}\,\mathfrak{H}\big)+\mu_0(\Nabla\hh)^\top\mm
-\frac\varrho2({\rm div}\,\vv)\,\vv+\varrho{\bm g}\,,
    \\\nonumber
    &\hspace*{4em}\text{where }\ 
    \bm{S}_\text{\sc e}^{}=\varphi_{\Ee}'(\Ee,\mm)+\psi(\Ee,\mm,\theta)\bbI
\,,\ \ \
\bm{S}_\text{\sc v}^{}=\DD\EE(\vv)\,,\ \ \hh=\hh_\text{\rm geo}{+}\hh_{\rm dem}\,,
    \\\nonumber
  &\hspace*{4em}\text{and }\ \ \ \bm{S}_\text{\sc c}^{}=
\!\!\!\lineunder{
\kappa\mu_0\Big(\Nabla\mm{\otimes}\Nabla\mm
-\frac1{\qexp}|\Nabla\mm|^2\bbI\Big)}{capillarity (Korteweg-like) stress}
\!\!\!\!\!-\!\!\!\linesunder{\mu_0{\rm skw}\big(\hh_{\rm drv}
{\otimes}\mm\big)}{couple stress due}{to magnetic dipoles}
\\&\hspace*{4em}\text{and }\ \ \ \hh_{\rm dem}=-\nabla u\,,\ \ \ 
\mathfrak{H}=\NU|\Nabla\EE(\vv)|^{p-2}\Nabla\EE(\vv)\,,
\label{MG-1}\\[-.1em]&
\ZJ\Ee+\ZJ\Ep=\EE(\vv)\,,
\label{MG-2}
\\\label{MG-3}
 &
 \GM(\theta)\ZJ\Ep +{\rm dev}\,\bm{S}_\text{\sc e}^{}=\varkappa\Delta\ZJ\Ep\,,
\\[-.3em]&\partial_{\ZJ\mm}\zeta(\theta;\ZJ\mm)
\ni\hh_{\rm drv}\ \ \text{ with }\ \ \hh_{\rm drv}=
\hh-\frac{\psi_\mm'(\Ee,\mm,\theta)}{\mu_0}+\kappa\Delta\mm\,,
\label{MG-4}
\\[-.2em]&\label{Max++}
\Delta u={\rm div}(\chi_\varOmega^{}\mm)\ \ \ \ \text{ on }\ \ \R^d\,,
\\[-.0em]&\nonumber
\pdt\W+{\rm div}\big(\W\vv{-}\KK(\theta)\Nabla\theta\big)=\xi\big(\theta;\EE(\vv),\ZJ\Ep,\ZJ{\mm}\big)
+\theta\CCC'(\mm){\cdot}\DT\mm+
\big(\theta\CCC(\mm){+}\phi(\theta)\big)\,{\rm div}\,\vv
\\[-.4em]&\nonumber
\hspace*{5em}\text{ with }\ \ \W=\gamma(\theta)
\ \ \text{ and }\ \ \xi\big(\theta;\EE,\ZJ\Ep,\ZJ{\mm}\big)
=\bm{S}_\text{\sc v}^{}{:}\EE+\mathfrak{H}\Vdots\Nabla\EE
+\GM(\theta)|\ZJ\Ep|^2
\\[-.2em]&\hspace*{21.9em}
+\varkappa|\Nabla\ZJ\Ep|^2+\mu_0\partial_{\ZJ\mm}\zeta(\theta;\ZJ\mm){\cdot}\ZJ\mm \,.
 \label{MG-6+thermo}
 \end{align}\end{subequations}
The equations (\ref{MG}a-d,f) are considered on the domain $\varOmega$, in
contrast to \eqref{Max++} which is considered on the whole $\R^d$.
Of course, the corotational derivatives $\ZJ\Ee$ and $\ZJ\Ep$ are from
\eqref{ZJ-tensor} and $\ZJ\mm$ is from \eqref{ZJ-vector}. Note that,
in view of \eqref{free}, we could equivalently  write
$\bm{S}_\text{\sc e}^{}=\psi_{\Ee}'(\Ee,\mm,\theta)+\psi(\Ee,\mm,\theta)\bbI$.
For the capillarity-like term and
the skew-symmetric stress see also \cite{BePGVa01DDFG,DeSPoG95ISIS}.
 The magnetizaton flow rule \eqref{MG-4} with the corotational derivative
 $\ZJ\mm$  see also \cite{BePGVa01DDFG,Maug76CTDF} where it is articulated
 that the magnetization is ``frozen'' in the deforming medium if $\ZJ\mm=0$ 
which then means that the magnetization is transported and rotates at
the same local rate as the deforming medium; this is the situation below
the blocking temperature $\theta_{\rm b}$ and when the total driving field
$\hh_{\rm drv}$ has small magnitude.

The first term in \eqref{free} allows for involving also magnetostrictive
effects in general. Its simplest (and isotropic) form without
magnetostrictive effects (which can surely be neglected in paleomagnetism
in rocks) is 
\begin{align}\label{psi-special}
\varphi(\Ee,\mm)=
\frac d2K_\text{\sc e}|{\rm sph}\,\Ee|^2\!+G_\text{\sc e}|{\rm dev}\,\Ee|^2\!
+b_0|\mm|^4\!-a_0\theta_{\rm c}|\mm|^2\ \ \ \text{ and }\ \ \
\CC(\mm)=a_0|\mm|^2
\end{align}
with $K_\text{\sc e}$ is the bulk elastic modulus, $G_\text{\sc e}$ the shear
elastic modulus, and with some $a_0>0$ and $b_0$ and with $\theta_{\rm c}$
denoting the Curie temperature. Neglecting the demagnetizing field, i.e.\
for $\hh_{\rm dem}=\mathbf0$, the minimum of such $\psi(\Ee,\cdot,\theta)$ with
respect to $\mm$ is attained at the orbit $|\mm|=m_{\rm s}$ with the
saturation magnetization
\begin{align}
m_{\rm s}=m_{\rm s}(\theta)=\begin{cases}
\sqrt{a_0(\theta_{\rm c}{-}\theta)/b_0}&\text{for }\ \theta<\theta_{\rm c}\,,\\
\qquad0&\text{for }\ \theta\ge\theta_{\rm c}\,;\end{cases}
\end{align}
note that $\psi(\Ee,\cdot,\theta)$ from \eqref{psi-special} becomes
convex with the single minimizer $\mm={\bm 0}$ for $\theta$ above
the Curie temperature $\theta_{\rm c}$, cf.\ Figure~\ref{fig1}.

The further ingredient for building the model is the (pseudo)potential
of dissipative forces as a convex functional of rates. In our case, we will
choose it temperature dependent as
\begin{align}\label{dissip-pot}
(\vv,\ZJ\Ep,\ZJ\mm)\mapsto\frac{\DD}2|\EE(\vv)|^2+\frac{\NU}p|\Nabla\EE(\vv)|^p
+\frac{\GM(\theta)}2|\ZJ\Ep|^2+\frac\varkappa2|\nabla\ZJ\Ep|^2+\zeta(\theta;\ZJ\mm)\,.
\end{align}
For notational simplicity, we consider the same Stokes-type (hyper)viscosity
in the deviatoric and the spherical parts by considering single parameters
$\DD$ and $\NU$.
The simplest dissipation potential $\zeta=\zeta(\theta;\ZJ\mm)$ which
gives the typical hysteretic response in the $\mm/\hh$-diagrammes
is non-differentiable at the magnetization rate $\ZJ\mm=\mathbf0$:
\begin{align}
\zeta(\theta;\ZJ\mm)=h_{\rm c}(\theta)|\ZJ\mm|\ \ \ \text{ with }\ \
h_{\rm c}(\theta)=\begin{cases}\text{high}&
\text{for $\theta$ below $\theta_{\rm b}$},\\
\text{low (or zero)}\!\!\!&\text{for $\theta$ above $\theta_{\rm b}$,
cf.\ Fig.\,\ref{fig1}}.\end{cases}
\label{hc-special}\end{align}
For temperature below blocking temperature $\theta_{\rm b}$, the
subdifferential of $\zeta(\theta;\cdot)$ is then multivalued at
$\ZJ\mm=\mathbf0$, as depicted in Figure~\ref{fig2}a.
Together with \eqref{psi-special}, this simplest scenario \eqref{psi-special}
gives the hysteretic loops as in Figure~\ref{fig2}b.
\begin{center}
\begin{my-picture}{35}{17}{fig2}
\hspace*{-2em}{\includegraphics[width=40.5em]{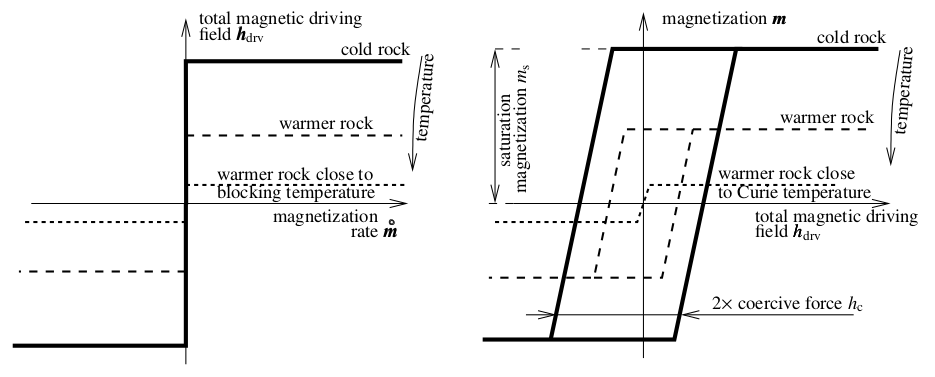}}
\end{my-picture}
\nopagebreak
\\
  {\small\sl\hspace*{-.5em}Fig.\,\ref{fig2}a:~\begin{minipage}[t]{14em}
Schematic multivalued function $\partial_{\ZJ\mm}\zeta(\theta;\cdot)$ as
a subdifferential of the potential $\zeta(\theta;\cdot)$ which is
nondifferentiable at $\ZJ\mm=\mathbf0$.
\end{minipage}}\hspace*{4em}{\small\sl\hspace*{-2em}Fig.\,\ref{fig2}b:~\begin{minipage}[t]{17em}
An idealized hysteresis in $\mm/\hh_\text{\rm geo}$-de\-pe\-ndence for various
temperatures during rate-independent thermoremanent magnetization (TRM).
\end{minipage}
} 
\end{center}
The slope of the hysteretic loops in Figure~\ref{fig2}b is
due to demagnetizing effects and depends on the shape of the magnets;
actually, the loops typically are rather curved than piecewise affine, 
cf.~e.g.\ \cite{KruRou03SSIH}. Actually, the model
\eqref{psi-special}--\eqref{hc-special} covers also 
the isothermal remanent magnetization (IRM) when the intensity of
the external magnetic field $\hh_\text{\rm geo}$ is sufficiently large. 

To cover also the viscous remanent magnetization (VRM), one should modify
the dissipation potential \eqref{hc-special} in a small neigbourhood of
$\mathbf0$. Conceptually, one can consider
\begin{align}
\zeta(\theta;\ZJ\mm)=h_{\rm c}(\theta)|\ZJ\mm|+\epsilon|\ZJ\mm|^\rexp+
\begin{cases}\tau_{\rm c}(\theta)|\ZJ\mm|^2&\text{for $|\ZJ\mm|\le m_{\rm r}$}\,,
\\\tau_{\rm c}(\theta)m_{\rm r}^2&\text{for $|\ZJ\mm|>m_{\rm r}$}
\end{cases}
\label{hc-special+}\end{align}
with $h_{\rm c}$ from \eqref{hc-special} and with some presumably large
$\tau_{\rm c}(\theta)$ and presumably very small $\epsilon>0$ and small
$m_{\rm r}>0$;
the physical unit of the magnetization rate $m_{\rm r}$ is Am$^{-1}$s$^{-1}$
while $\tau_{\rm c}$ is some time constant (in the dimension seconds) similarly
as $\epsilon$ if the exponent $\rexp$ would equal 2. 
 the potential $\zeta(\theta;\cdot)$ is strictly convex,
uniformly with respect to $\theta$, which will simplify the analysis.
A subdifferential of such a potential is illustrated in Figure~\ref{fig3}.
\begin{center}
\begin{my-picture}{35}{15}{fig3}
\hspace*{-0em}{\includegraphics[width=38em]{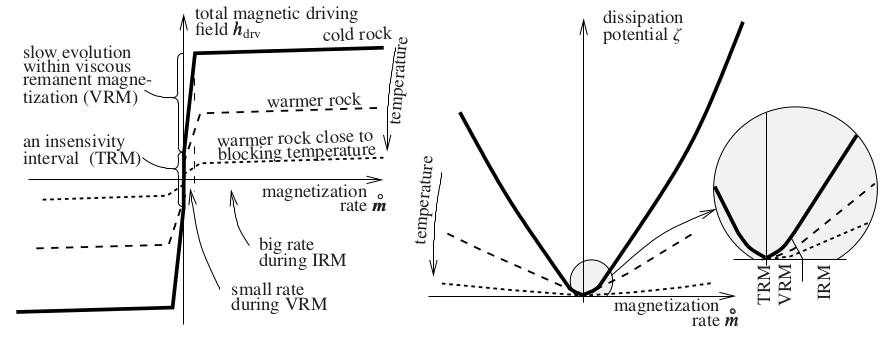}}
\end{my-picture}
\nopagebreak
\\
{\small\sl\hspace*{-.5em}Fig.\,\ref{fig3}a:~\begin{minipage}[t]{16.5em}
A schematic refinement of Fig.\,\ref{fig2}a involving also
viscous remanent magnetization (VRM) together with thermoremanent
magnetization (TRM) and rate-dependent isothermal
remanent magnetization (IRM).
\end{minipage}}\hspace*{3em}
{\small\sl\hspace*{-2em}Fig.\,\ref{fig3}b:~\begin{minipage}[t]{15em}
The corresponding dissipation potential
non-differentiable at magnetization rate zero and depending on temperature.
\end{minipage}
} 
\end{center}

An important modelling assumption that the elastic-strain tensor $\Ee$ is
kept
symmetric during the evolution, which is granted by initial conditions
on both $\Ee$ and $\Ep$ and by the choice of the objective time derivative
as corotational, cf.\ the first property in \eqref{ZJ-property}. Then also
the elastic-stress tensor $\varphi_\Ee'$ is symmetric. Moreover, we will
use the second property in \eqref{ZJ-property} and the initial condition
on $\Ep$ to keep $\Ep$ not only symmetric but also trace-free, i.e.\
deviatoric.

\def\bbM{\mathbb M}
\begin{remark}[{\sl Jeffreys and Kelvin-Voigt rheologies}]\upshape
To reveal the rheological model behind (\ref{MG}a--c), assuming for a moment
that $\varkappa=0$ and $\NU=0$, we can eliminate the internal
variables $\Ee$ and $\Ep$ and obtain a relation between the total stress
$\bm{S}$ and total strain rate $\EE(\vv)$. Denoting by $\bbC$, $\bbD$,
and $\bbM$ the tensors of elasticity modulus, the Stokes-viscosity modulus,
and the Maxwellian viscosity modulus, respectively, we can write
$\bm{S}_\text{\sc e}^{}=\bbC\Ee=\bbD\ZJ\Ep$ and
$\bm{S}_\text{\sc v}^{}=\bbD\EE(\vv)$. Using also \eqref{MG-2}, by some algebraic
manipulation, we obtain
\begin{align}\nonumber\\[-2.2em]\label{Jeffreys}
\bbC^{-1}\ZJ{\bm S}+\bbM^{-1}{\bm S}
=\bbC^{-1}\bbD\ZJ{\overline{\EE(\vv)}}+(\bbI{+}\bbM^{-1}\bbD)\EE(\vv)\,;
\end{align}
at small strains cf.\ \cite[Sect.\,6.6]{KruRou19MMCM}. This is a relation for
the Jeffreys rheology written in Eulerian coordinates. For vanishing $\bbD$,
\eqref{Jeffreys} degenerates to
\begin{align}\nonumber\\[-2.2em]\label{Maxwell}
\bbC^{-1}\ZJ{\bm S}+\bbM^{-1}{\bm S}=\EE(\vv)\,,
\end{align}
which is the relation for the {\it Maxwell rheology}. For $\bbM$ very large
``going to infinity'', i.e.\ for $\bbM^{-1}$ vanishing, \eqref{Jeffreys}
degenerates to
\begin{align}\nonumber\\[-2.6em]
\label{KV}
\ZJ{\bm S}=\bbD\ZJ{\overline{\EE(\vv)}}+\bbC\EE(\vv)\,,
\end{align}
which is the relation for the {\it Kelvin-Voigt rheology}.
In terms of acceleration $\DT\vv$, one can also express
$\ZJ{\overline{\EE(\vv)}}$ occurring in \eqref{Jeffreys} and \eqref{KV} as 
$\EE(\DT\vv)-{\rm sym}(\Nabla\vv\Nabla\vv)+{\rm skw}(\Nabla\vv^\top\Nabla\vv)$.
In \eqref{psi-special} and \eqref{dissip-pot}, the above moduli can be written
componentwise as
$[\bbC]_{ijkl}=K_\text{\sc e}\delta_{ij}\delta_{kl}+G_\text{\sc e}d_{ijkl}$
with $d_{ijkl}=\delta_{ik}\delta_{jl}{+}
\delta_{il}\delta_{jk}{-}2\delta_{ij}\delta_{kl}/d$ with $\delta$ the Kronecker
symbol and with $K_\text{\sc e}$ and $G_\text{\sc e}$ from \eqref{psi-special},
$[\bbD]_{ijkl}=\DD d_{ijkl}$, and $\bbM=\GM(\theta)d_{ijkl}$. This means that our
model combines temperature-dependent Jeffreys rheology in the deviatoric part
with the Kelvin-Voigt rheology in the spherical (volumetric) part.
\end{remark}

\begin{remark}[{\sl Ellimination of $\Ep$}]\label{rem-R}\upshape
In fact, $\Ep$ itself occurs in \eqref{MG} only as its corotation rate
$\ZJ\Ep$. Thus, introducing the new variable $\RR:=\ZJ\Ep$, one can write
\eqref{MG} in terms of $\RR$ without $\Ep$. In particular, the parabolic
equation \eqref{MG-3} turns into the elliptic equation
$\GM(\theta)\RR+{\rm dev}\,\bm{S}_\text{\sc e}^{}=\varkappa\Delta\RR$
and one does not need any initial condition for $\Ep$. Knowing a solution
$(\vv,\RR,\mm,u,\theta)$ of such system and prescibing an initial
condition for $\Ep$, one can then reconstruct the inelastic strain $\Ep$
by solving the nonhomogeneous transport equation $\DT\Ep=\RR+{\rm skw}(\nabla\vv)\Ep-\Ep{\rm skw}(\nabla\vv)$.
\end{remark}

\section{Energetics and thermodynamics of the model}\label{sec-thermodyn}

To reveal the energetics behind the system, we must specify some boundary
conditions for the particular equations (\ref{MG}a,c--f), say
\begin{subequations}\label{BC}
  \begin{align}\label{BC1}
&\vv{\cdot}\nn=0\,,\ \ \ \ \Nabla\EE(\vv){:}(\nn{\otimes}\nn)=\mathbf{0}\,,
\\[-.0em]\label{BC1-}
&\big[\big({\bm S}_\text{\sc e}^{}\!+\bm{S}_\text{\sc v}^{}\!
    +\bm{S}_\text{\sc c}^{}-\mu_0(\hh{\cdot}\mm)\bbI
  -{\rm div}\mathfrak{H}\big)\nn+\divS(\mathfrak{H}\nn)
  \big]_\text{\sc t}^{}\!=\mathbf{0}\,,\ \ \
    \\[-.2em]\label{BC2}
   &
  (\nn{\cdot}\Nabla)\ZJ\Ep=\mathbf{0}\,,\ \ \ \
 (\nn{\cdot}\Nabla)\mm=\mathbf{0}\,,\ \ \ \ 
 \lim_{|x|\to\infty}u(x)=0\,,\ \ \text{ and }\ \
 \\[-.3em]
 &\label{BC5}\nn{\cdot}\KK(\theta)\Nabla\theta=j_{\rm ext}\,,
\end{align}
\end{subequations}
where $\divS={\rm tr}(\NablaS)$ with ${\rm tr}(\cdot)$ being the trace of a
$(d{-}1){\times}(d{-}1)$-matrix, denotes the $(d{-}1)$-dimensional
surface divergence and $\NablaS v=\Nabla v-\frac{\partial v}{\partial\nn}\nn$ 
being the surface gradient of $v$.
The first condition in \eqref{BC1} allows for considering a fixed
domain $\varOmega$, otherwise the situation would be extremely complicated
because, beside other analytical complications, evolving domains might exhibit
selfpenetration. The condition in \eqref{BC1-} employs the notation
$[\,\cdot\,]_\text{\sc t}^{}$ for the tangential component of a vector.
The last condition involves an external heat flux $j_{\rm ext}$ considered
prescribed.

To see the  magneto-mechanical energy balance behind the system
(\ref{MG}a-e) is seen by testing the equation
\eqref{MG-1} by $\vv$, \eqref{MG-3} by $\ZJ\Ep$, 
\eqref{MG-4} by $\mu_0\ZJ\mm$, and \eqref{Max++} by $\mu_0\pdt{}u$.
In this section, we perform the calculus only formally, i.e.\ assuming
that the system \eqref{MG} with the boundary conditions \eqref{BC}
has a solution which is sufficiently smooth.

Executing the first mentioned test of \eqref{MG-1} by $\vv$,
the inertial terms integrated over $\varOmega$ gives, by using the Green
formula and the boundary condition $\vv{\cdot}\nn=0$, that 
\begin{align}\nonumber
&\int_\varOmega\varrho\Big(\DT\vv+\frac12\vv\,{\rm div}\,\vv\Big)
{\cdot}\vv\,\d x
=\int_\varOmega\varrho\Big(\pdt\vv+(\vv{\cdot}\Nabla)\vv
+\frac12\vv\,{\rm div}\,\vv\Big){\cdot}\vv\,\d x
\\[-.1em]\nonumber&\ 
=\frac{\d}{\d t}\int_\varOmega\frac\varrho2|\vv|^2\,\d x
+\!\int_\varOmega\frac\varrho2|\vv|^2{\rm div}\,\vv
-\varrho\vv{\cdot}{\rm div}(\vv{\otimes}\vv)-|\vv|^2(\vv{\cdot}\!\!\!\!\lineunder{\!\!\Nabla\varrho}{$=\mathbf{0}$}\!\!\!\!)\,\d x
+\int_\varGamma\frac\varrho2|\vv|^2(\!\!\!\!\lineunder{\vv{\cdot}\nn_{_{}}}{=0}\!\!\!\!)\,\d S
\\[-.3em]&\
=\frac{\d}{\d t}\int_\varOmega\frac\varrho2|\vv|^2\,\d x
+\int_\varOmega\varrho\Big(\frac12|\vv|^2{\rm div}\,\vv
-\vv{\cdot}{\rm div}(\vv{\otimes}\vv)\Big)\,\d x
=\frac{\d}{\d t}\int_\varOmega\frac\varrho2|\vv|^2\,\d x\,,
\label{formula1}\end{align}
where we used also that $\int_\varOmega\vv{\cdot}{\rm div}(\vv{\otimes}\vv)\,\d x
=\int_\varOmega|\vv|^2{\rm div}\,\vv+(\vv{\otimes}\vv){:}\Nabla\vv\,\d x
=\int_\varGamma|\vv|^2(\vv{\cdot}\nn)\,\d S+\int_\varOmega|\vv|^2{\rm div}\,\vv
-\vv{\cdot}{\rm div}(\vv{\otimes}\vv)\,\d x
=\int_\varOmega\frac12|\vv|^2{\rm div}\,\vv\,\d x$ by the Green formula with the
boundary condition $\vv{\cdot}\nn=0$.
For the stress terms in \eqref{MG-1}, we use simply the Green theorem.
The peculiarity is behind the 2nd-grade hyperstress for which we should
use twice the Green formula over $\varOmega$ and once a surface Green formula
over the boundary $\varGamma$, cf.\ \cite[Sect.\,2.4.4]{Roub13NPDE} for
details. The other stress terms need only once Green theorem, which leads to
$\int_\varOmega\bm{S}_\text{\sc e}^{}{:}\EE(\vv)\,\d x$,
$\int_\varOmega\bm{S}_\text{\sc v}^{}{:}\EE(\vv)\,\d x$, and
$\int_\varOmega{\bm S}_\text{\sc c}^{}{:}\Nabla\vv\,\d x$. The first term
$\bm{S}_\text{\sc e}^{}{:}\EE(\vv)$ is balanced with terms resulting from
\eqref{formula2} and \eqref{formula4} below, while the second term
$\bm{S}_\text{\sc v}^{}{:}\EE(\vv)$
contributes to the dissipation rate $\xi$ in \eqref{MG-6+thermo}.
The last term ${\bm S}_\text{\sc c}^{}{:}\Nabla\vv$ is 
balanced with the particular terms arising in 
\eqref{formula3} and \eqref{calculus-PM4}.

The flow rule \eqref{MG-3} for $\Ep$ is to be multiplied by $\ZJ\Ep$.
From the term ${\rm dev}\bm{S}_\text{\sc e}^{}$, we obtain
\begin{align}\nonumber
  &{\rm dev}\bm{S}_\text{\sc e}^{}{:}\ZJ\Ep=\big(\bm{S}_\text{\sc e}^{}
  -{\rm sph}\bm{S}_\text{\sc e}^{}\big){:}\ZJ\Ep
=\bm{S}_\text{\sc e}^{}{:}\ZJ\Ep=\bm{S}_\text{\sc e}^{}{:}\big(\EE(\vv)
-\ZJ\Ee\big)
\\&\qquad\qquad\nonumber
=\bm{S}_\text{\sc e}^{}{:}\EE(\vv)-
\bm{S}_\text{\sc e}^{}{:}\Big(\pdt\Ee+(\vv{\cdot}\Nabla)\Ee\Big)
-\bm{S}_\text{\sc e}^{}{:}\big(\Ee\,{\rm skw}(\Nabla\vv)
-{\rm skw}(\Nabla\vv)\Ee\big)
\\&\qquad\qquad\nonumber
=\bm{S}_\text{\sc e}^{}{:}\,\EE(\vv)
-\bm{S}_\text{\sc e}^{}{:}\Big(\pdt\Ee+(\vv{\cdot}\Nabla)\Ee\Big)
+{\rm skw}\big(\hspace*{-1.7em}\linesunder{\bm{S}_\text{\sc e}^{}\Ee^\top\!\!
-\Ee^\top\bm{S}_\text{\sc e}}{=\,$\mathbf0$ as $\Ee$ and $\bm{S}_\text{\sc e}$ are\ \ }
{symmetric and commute}\hspace*{-1.7em}\big){:}\Nabla\vv
\\[-.1em]&\qquad\qquad
=\bm{S}_\text{\sc e}^{}
{:}\EE(\vv)-\Big(\varphi_\Ee'(\Ee,\mm)+\psi(\Ee,\mm,\theta)\bbI\Big){:}\Big(\pdt{\Ee\!}+(\vv{\cdot}\Nabla)\Ee\Big)\,.
\label{formula2}\end{align}
The second equality
in \eqref{formula2} has used the orthogonality of the deviatoric and
the spherical parts of $\bm{S}_\text{\sc e}^{}$ and the mentioned
attribute that the inelastic strain $\Ep$ is deviatoric during the whole
evolution. The third equality has used \eqref{MG-2}. The penultimate equality
has used the matrix algebra $A{:}(BC)=(B^\top\!A){:}C=(AC^\top){:}B$ while for
the last equality in \eqref{formula2} we used
$\bm{S}_\text{\sc e}^{}=\varphi_{\Ee}'(\Ee,\mm){+}\psi(\Ee,\mm,\theta)\bbI$
together with the mentioned symmetry of both $\Ee$ and $\bm{S}_\text{\sc e}^{}$
and the isotropy of $\varphi$ so that $\bm{S}_\text{\sc e}^{}$ and
$\Ee$ commute with each other, and as a result the skew-symmetric part
vanishes.

The test of \eqref{MG-4} by $\mu_0\ZJ\mm$ is quite technical. The term
$\kappa\Delta\mm{\cdot}\mu_0(\vv{\cdot}\Nabla)\mm$ 
is to be handled by using Green's formula twice. Namely,
\begin{align}\nonumber
\int_\varOmega\!\kappa\Delta\mm{\cdot}\mu_0(\vv{\cdot}\Nabla)\mm\,\d x
&=\!\int_\varGamma\kappa
(\nn{\cdot}\Nabla)\mm{\cdot}\mu_0\big((\vv{\cdot}\Nabla)\mm\big)\,\d S
\\[-.4em]&\nonumber\hspace*{2em}
-\!\int_\varOmega\mu_0\kappa\Nabla^2\mm{:}(\vv{\otimes}\Nabla\mm)+
\mu_0\kappa(\Nabla\mm{\otimes}\Nabla\mm){:}\EE(\vv)\,\d x
\\&\nonumber
=\int_\varGamma
\mu_0\kappa\Big((\nn{\cdot}\Nabla)\mm{\cdot}\big((\vv{\cdot}\Nabla)\mm\big)-
\frac12|\Nabla\mm|^2\vv{\cdot}\nn\Big)\,\d S
\\[-.4em]&\hspace*{2em}+\int_\varOmega
\frac{\mu_0\kappa}{\qexp}|\Nabla\mm|^{\qexp}{\rm div}\,\vv
-\mu_0\kappa(\Nabla\mm{\otimes}\Nabla\mm){:}\EE(\vv)\,\d x
\,,
\label{test-damage}\end{align}
where the boundary integral vanishes due to the boundary conditions
$(\nn{\cdot}\Nabla)\mm=\mathbf0$ and $\vv{\cdot}\nn=0$.
In the last integral, we see the terms which (counting the factor
$\mu_0\kappa$) are balanced by the Korteweg-like stress in \eqref{MG-1}.
Moreover, we use the Green theorem also for the
driving magnetic field $\hh=\hh_\text{\rm geo}{+}\hh_\text{\rm dem}$:
\begin{align}
\int_\varOmega\hh{\cdot}(\vv{\cdot}\Nabla)\mm\,\d x
=\int_\varGamma(\hh{\cdot}\mm)\vv{\cdot}\nn\,\d S-
\int_\varOmega(\Nabla\hh)^\top\mm{\cdot}\vv
+(\hh{\cdot}\mm)\,{\rm div}\,\vv\,\d x\,,
\label{formula2+}\end{align}
which (counting the factor $\mu_0$) is balanced with the force
$\mu_0(\Nabla\hh)^\top\mm$ and the pressure $\mu_0\hh{\cdot}\mm$ in
\eqref{MG-1}. This is used to handle the terms in 
$\mu_0\hh_{\rm drv}{\cdot}\ZJ\mm$
by the calculus
\begin{align}\nonumber
&\int_\varOmega\!\mu_0\hh_{\rm drv}{\cdot}\ZJ\mm\,\d x
=\nonumber\!\int_\varOmega\!\big(\mu_0 \kappa\Delta\mm
-\psi_{\mm}'(\Ee,\mm,\theta)+\mu_0(\hh_{\rm geo}{+}\hh_{\rm dem})\big){\cdot}\DT\mm
-\mu_0\hh_{\rm drv}{\cdot}{\rm skw}(\Nabla\vv){\mm}\,\d x
  \\[-.1em]&\qquad\ =\nonumber
\frac{\d}{\d t}\int_\varOmega\!\mu_0\hh_\text{\rm geo}{\cdot}\mm
-\frac{\mu_0\kappa}{\qexp}|\Nabla\mm|^{\qexp}\,\d x
\\[-.3em]&\nonumber\hspace*{4.7em}
+\!\int_{\varGamma}\mu_0\kappa\DT\mm{\cdot}
\!\!\!\!\lineunder{(\nn{\cdot}\Nabla)\mm}{$=\mathbf0$}\!\!\!
-\frac{\mu_0\kappa}\qexp|\Nabla\mm|^{\qexp}\!\!\!\!\lineunder{\vv{\cdot}\nn}{$=0$}\!\!\!\!
+\mu_0(\hh_{\rm drv}{\cdot}\mm)(\!\!\!\!\lineunder{\vv{\cdot}\nn}{$=0$}\!\!\!\!)\,\d S
\\[-.7em]&\nonumber\hspace*{4.7em}
+\!\int_\varOmega\!\bigg(\frac{\mu_0\kappa}{\qexp}|\Nabla\mm|^{\qexp}{\rm div}\,\vv
-\mu_0\kappa
(\Nabla\mm{\otimes}\Nabla\mm){:}\EE(\vv)
-\varphi_{\mm}'(\Ee,\mm){\cdot}\frac{\partial\mm}{\partial t}
\\[-.3em]&\nonumber\hspace*{6.7em}
+\big(\mu_0\hh_{\rm dem}{-}\theta\CCC'(\mm)\big){\cdot}\DT\mm
-\mu_0\pdt{\hh_\text{\rm geo}}{\cdot}\mm
-\varphi_{\mm}'(\Ee,\mm){\cdot}(\vv{{\cdot}}\Nabla)\mm
\\[-.4em]&\hspace*{6.7em}
-\mu_0(\nabla\hh)^{\top}\mm{\cdot}\vv-\mu_0(\hh{\cdot}\mm){\rm div}\,\vv
-\mu_0{\rm skw}\big(\hh_{\rm drv}{\otimes}\mm\big){:}\Nabla\vv\!\bigg)\,\d x\,,
\label{formula3}
\end{align}
Here we see the skew-symmetric contribution
${\rm skw}(\hh_{\rm drv}{\otimes}\mm)$
which is balanced by the skew-symmetric couple-like stress in \eqref{MG-1},
as well as the pressure-like term $\hh{\cdot}\mm$ and the force
$(\Nabla\hh)^\top\mm$ in \eqref{MG-1}.
The term $\theta\CCC'(\mm){\cdot}\DT\mm$ will occur in \eqref{energy+}
and will be balanced in the total energy with the adiabatic heat source/sink 
in  \eqref{MG-6+thermo}. 

The partial time derivative terms in \eqref{formula2} and \eqref{formula3}
can be merged by the calculus as
\begin{align}
&\varphi_\Ee'(\Ee,\mm){:}\pdt\Ee+\varphi_\mm'(\Ee,\mm){:}\pdt\mm=
\pdt{}{\varphi(\Ee,\mm)}\,.
\label{formula4-}\end{align}
The convective terms in \eqref{formula2} and \eqref{formula3} can be merged
for the calculus and by the Green formula
\begin{align}\nonumber
&\int_\varOmega\varphi_{\Ee}'(\Ee,\mm):(\vv{\cdot}\Nabla)\Ee
  +\varphi_{\mm}'(\Ee,\mm){\cdot}(\vv{\cdot}\Nabla)\mm\,\d x
\\[-.4em]&\qquad
=\int_\varOmega\Nabla\varphi(\Ee,\mm){\cdot}\vv\,\d x
 =\int_\varGamma
 \varphi(\Ee,\mm)\!\!\!\lineunder{\vv{\cdot}\nn}{=0}\!\!\!\!\!\,\d S-
\int_\varOmega\varphi(\Ee,\mm){\rm div}\,\vv\,\d x\,.
\label{formula4}\end{align}
Here we can see the term $\varphi(\Ee,\mm)$ which are balanced by the
pressure contribution in the elastic stress $\bm{S}_\text{\sc e}^{}$.

For \eqref{Max++} tested by $\mu_0\pdt{}u$, we use the calculus,
including the Green theorem for the convective term, to obtain

\begin{align}\nonumber
0&=\int_{\R^d}\mu_0{\rm div}\big(\chi_\varOmega^{}\mm-\Nabla u\big)\pdt{u}\,\d x=
\int_{\R^d}\mu_0\big(\Nabla u-\chi_\varOmega^{}\mm\big){\cdot}\Nabla\pdt{u}\,\d x
\\\nonumber
&=\frac{\d}{\d t}\int_{\R^d}\frac{\mu_0}2|\Nabla u|^2\,\d x
+\int_\varOmega\mu_0\mm{{\cdot}}\pdt{\hh_{\rm dem}}\,\d x
\\\nonumber&=\frac{\d}{\d t}\bigg(\int_{\R^d}\frac{\mu_0}2|\Nabla u|^2
+\int_\varOmega\mu_0\mm{\cdot}\hh_{\rm dem}\,\d x\bigg)
-\int_\varOmega\mu_0\pdt\mm{\cdot}\hh_{\rm dem}\,\d x
\\\nonumber
&=\frac{\d}{\d t}\bigg(\int_{\R^d}\!\!\frac{\mu_0}2|\Nabla u|^2\,\d x
+\!\int_\varOmega\!\mu_0\mm{\cdot}\hh_{\rm dem}\,\d x\!\bigg)
+\!\!\int_\varOmega\!\!
\mu_0\Big({\rm skw}(\Nabla\vv)\mm-\ZJ\mm
-(\vv{\cdot}\Nabla)\mm\Big){\cdot}\hh_{\rm dem}\,\d x
\\\nonumber&=\frac{\d}{\d t}\bigg(\int_{\R^d}\!\!\frac{\mu_0}2|\Nabla u|^2\,\d x
+\!\int_\varOmega\mu_0\mm{\cdot}\hh_{\rm dem}\,\d x\!\bigg)
-\int_\varGamma\mu_0(\hh_{\rm dem}{\cdot}\mm)\!\!\!\lineunder{\vv{\cdot}\nn}{$=0$}\!\!\!\d S
\\[-.4em]\nonumber&\hspace{3em}
+\int_\varOmega\Big(\mu_0{\rm skw}(\hh_{\rm dem}{\otimes}\mm){:}\Nabla\vv
-\mu_0\ZJ\mm{\cdot}\hh_{\rm dem}
\\[-.6em]&\hspace{16em}
+\mu_0\vv{\cdot}(\Nabla\hh_{\rm dem})^\top\mm
+\mu_0(\hh_{\rm dem}{\cdot}\mm){\rm div}\,\vv\Big)\,\d x\,.
\label{calculus-PM4}\end{align}
The last integral reveals how the force $\mu_0(\Nabla\hh_{\rm dem})^\top\mm$ and
the stress contribution
$\mu_0{\rm skw}(\hh_{\rm dem}{\otimes}\mm)-\mu_0(\hh_{\rm dem}{\cdot}\mm)\bbI$
arises energetically in the momentum equation.
The term $\mu_0\ZJ\mm{\cdot}\hh_{\rm dem}$ is canceled with the term
$\mu_0\Nabla u\cdot\ZJ\mm$ arising from \eqref{MG-4} tested by $\mu_0\ZJ\mm$.

To summarize the above calculations, we state:

\begin{proposition}[{\sl Magneto-mechanical energy balance}] Any
smooth solution of the system (\ref{MG}a--e) with the boundary conditions \eqref{BC}
satisfies the identity
\begin{align}\nonumber
  &\frac{\d}{\d t}
  \int_\varOmega\!\!\!\!\linesunder{\frac\varrho2|\vv|^2}{kinetic}{energy}\!\!\!\!+\!\!\!\!\linesunder{\varphi(\Ee,\mm)+
    \frac{\kappa\mu_0}{\qexp}|\Nabla\mm|^{\qexp}}{stored}{energy}\!\!\!\!
    +\!\!\!\!\linesunder{\mu_0\hh_\text{\rm geo}{\cdot}\mm_{_{_{_{_{_{}}}}}}\!\!\!}{Zeeman}{energy}\!\!\,\d x+\int_{\R^d}\!\!\!\!\!\!\!\!\!\!\linesunder{\frac{\mu_0}2|\Nabla u|^2_{_{_{_{_{_{}}}}}}\!\!\!}{energy of de-}{magnetizing field}\!\!\!\!\!\!\!\,\d x
\\&\nonumber\hspace*{0em}
+\int_\varOmega\!\!\!\!\linesunder{\xi\big(\theta;\EE(\vv),\ZJ\Ep,\ZJ{\mm}\big)
    _{_{_{_{}}}}\!\!}{dissipation rate}{from \eqref{MG-6+thermo}}\!\!\!
    \d x=\int_\varOmega\!\!\!\!\linesunder{\varrho\ff{\cdot}\vv
+\frac{\!\partial \hh_\text{\rm geo}\!}{\partial t}{\cdot}\mm
  }{power of}{external load}\!\!\!\!-\!\!\!\!\!\lineunder{\theta\CCC'(\mm){\cdot}\DT\mm
    -\big(\theta\CCC(\mm){+}\phi(\theta)\big){\rm div}\,\vv_{_{_{_{_{_{}}}}}}\!}{adiabatic effects}\!\!\!\!\d x\,.
\\[-1.7em]&  \label{energy+}
\end{align}
\end{proposition}

Adding \eqref{MG-6+thermo} tested by 1 and using also the boundary condition in
\eqref{BC5}, the adiabatic effects are absorbed in the internal heat energy
$\W$ and we obtain:

\begin{proposition}[{\sl Total energy balance}] Any
smooth solution of the system \eqref{MG} with the boundary conditions
\eqref{BC} satisfies the identity
\begin{align}\nonumber
  &\frac{\d}{\d t}\bigg(
  \int_\varOmega\!\!\!\!\linesunder{\frac\varrho2|\vv|^2}{kinetic}{energy}\!\!\!\!\!+\!\!\!\!\!\linesunder{\varphi(\Ee,\mm)+
  \frac{\kappa\mu_0}{\qexp}|\Nabla\mm|^{\qexp}
   }{stored}{energy}\!\!\!\!
+\!\!\!\!\linesunder{\mu_0\hh_\text{\rm geo}{\cdot}\mm_{_{_{_{_{_{}}}}}}\!\!\!}{Zeeman}{energy}\!\!
+\!\!\!\!\!\!\!\!\!\!\linesunder{\W_{_{_{_{_{}}}}}\!\!\!}{internal}{heat energy}
\!\!\!\!\!\!\!\!\!\,\d x
\\&\hspace*{6em}
+\int_{\R^d}\!\!\!\!\!\!\!\!\!\!\linesunder{\frac{\mu_0}2|\Nabla u|^2_{_{_{_{_{_{}}}}}}\!\!\!}{energy of de-}{magnetizing field}\!\!\!\!\!\!\!\,\d x\bigg)
=\int_\varOmega\!\!\!\!\linesunder{\varrho\ff{\cdot}\vv
  +\frac{\partial \hh_\text{\rm geo}\!}{\partial t}{\cdot}\mm}{power of
  external}{bulk load}\!\!\!\d x
+\int_\varGamma\!\!\!\!\!\!\!\!\!\!\!\!\!\!\!\linesunder{{j_{\rm ext}}_{_{_{_{_{_{}}}}}}\!\!\!}{power of}{external heating}\!\!\!\!\!\!\!\!\!\!\!\!\!\d S\,.
\label{energy-mag}\end{align}
\end{proposition}

The thermodynamical context of this model relies on an additive splitting
of the specific free energy $\psi$ into the purely mechanical part and the
thermal part $\phi$, cf.\ \eqref{psi-special}.
An important attribute of the model, beside keeping the energetics
\eqref{energy+}--\eqref{energy-mag}, is the entropy imbalance, i.e.\ 
the Clausius-Duhem inequality. The specific entropy is then
$\eta=-\psi_\theta'(\mm,\theta)=\phi'(\theta)-\CCC(\mm)$ with $\CCC$,
being  an extensive variable (in JK$^{-1}$m$^{-3}$) and transported by the
{\it entropy equation}
\begin{align}\pdt\eta+{\rm div}\big(\vv\,\eta\big)
=\frac{\xi-{\rm div}\,{\bm j}}\theta\qquad\qquad\qquad
\label{entropy-eq}\end{align}
with $\xi$ denoting the heat production rate from \eqref{MG-6+thermo}
and ${\bm j}$ the heat flux (here governed by the Fourier law
${\bm j}=-\KK\Nabla\theta$). 
Substituting $\eta=\phi'(\theta)-\CCC(\mm)$ into \eqref{entropy-eq}, we obtain
the heat-transfer equation 
\begin{align}
c(\theta)\DT\theta=\xi-{\rm div}\,{\bm j}
+\theta\big(\CCC'(\mm){\cdot}\DT\mm-\eta\,{\rm div}\,\vv\big)\ \ \
\text{ with the heat capacity }\
c(\theta)=\theta\phi''(\theta)\,;
\label{heat-eq-}\end{align}
note that temperature (in K) is an intensive variable and is transported
by the material derivative. There are the adiabatic heat source/sink 
terms $\theta\CCC'(\mm){\cdot}\DT\mm$ and 
$\theta\eta\,{\rm div}\,\vv=\theta(\phi'(\theta){-}
\CCC(\mm)){\rm div}\,\vv$ on the right-hand side due to
both the magnetic phase transition and the compressibility of the continuum.
Note that only the convective time derivative $\DT\mm$ but not the rotation
${\rm skw}(\Nabla\vv)\mm$ influences the adiabatic magnetic heat.

Furthermore, the {\it internal energy}
is given by the Gibbs relation $e=\psi+\theta\eta$, and splits here into the
purely magneto-mechanical part and the purely thermal part
$\W=\theta\phi'(\theta)-\phi(\theta)=:\gamma(\theta)$, namely:
\begin{align}
e=\!\!\!\!\lineunder{\varphi(\Ee,\mm)+\theta\CC(\mm)
-\phi(\theta)}{$=\psi(\Ee,\mm,\theta)$}
\!\!\!\!+\!\!\!\!
\lineunder{\theta\phi'(\theta)-\theta
\CCC(\mm)}{$=\theta\eta(\mm,\theta)$}
\!\!\!=\varphi(\Ee,\mm)+\W\,.
\end{align}
Since $\gamma'(\theta)=\theta\phi''(\theta)$, the heat equation
\eqref{heat-eq-} can be written as
\begin{align}
\DT\W=\xi-{\rm div}\,{\bm j}
+\theta\CCC'(\mm){\cdot}\DT\mm
+\theta\big(\CCC(\mm){-}\phi'(\theta)\big)\,{\rm div}\,\vv\,.
\label{heat-eq--}\end{align}
Alternatively, we can write  \eqref{heat-eq--} as 
\begin{align}
\pdt\W+{\rm div}\big(\W\vv\big)=\xi-{\rm div}\,{\bm j}
+\theta\CCC'(\mm){\cdot}\DT\mm
  +\big(\theta\CCC(\mm){+}\phi(\theta)\big){\rm div}\vv\,,
\label{heat-eq}\end{align}
which reflects the property that the thermal internal energy $\W$ in Jm$^{-3}$
is again an extensive variable and is transported like \eqref{entropy-eq}.
This reveals the structure of \eqref{MG-6+thermo}. In particular, from
\eqref{entropy-eq} we can see at least formally by the usual calculus,
relying on positivity of temperature:

\begin{proposition}[{\sl Clausius-Duhem entropy inequality}]
Assuming $\xi\ge0$ and $\KK\ge0$, any smooth solution of the system
(\ref{MG}a--e) with positive temperature $\theta$ with the boundary
conditions \eqref{BC} satisfies the identity
\begin{align}\nonumber
\frac{\d}{\d t}\int_\varOmega\eta\,\d x\ge
\int_\varGamma\!\!\!\!\lineunder{(\vv\eta-{\bm j}/\theta)}{entropy flux}\!\!\!\!\cdot\nn\,\d S=-\!\int_\varGamma\frac{j_{\rm ext}}\theta\,\d S\,.
\end{align}
\end{proposition}
The mentioned positivity of temperature can be ensured by non-negativity
of the boundary hear flux $j_{\rm ext}$ and positivity of the initial temperature
together with at least linear decay of $\phi(\cdot)$ to zero for $\theta\to0+$.

\begin{remark}[{\sl Oberbeck-Bousinesq buoyancy enhancement}]\label{OB}\upshape
An important phenomenon is that hot or even molten rocks (in particular magma)
are lighter than solid rocks, which rise a tendency of magma floating up in
the gravitation field $\ff$.
This can be included in a simplified way by the Oberbeck-Bousinesq
buoyancy model, usually used for incompressible media
which, anyhow, exhibit a slight thermal
  expansibility. This gives rise to an extra force by replacing
  $\ff$ with $\ff(1{-}b(\theta))$ in \eqref{MG-1} 
  with some $b(\cdot)$ continuous. Then \eqref{MG-6+thermo} expands 
  by the adiabatic heat source/sink term $b(\theta)\vv{\cdot}\ff$
and the analysis in Sect.\,\ref{sec-anal} below can easily be enhanced like
(even more  easily) for the term
$(\theta\CCC(\mm){+}\phi(\theta)){\rm div}\,\vv$ in \eqref{MG-6+thermo}.
\end{remark}

\begin{remark}[{\sl More general free energies}]\label{gen-free}\upshape
  The free energy \eqref{free} has a specific feature that $\psi_{\theta\theta}''$
  is independent of $\mm$. This leads to a specific ``transition''
  temperature (Curie or
  N\'eel) at which the spontaneous (saturation) magnetization
  $m_{\rm s}$ sharply falls to zero, see Figure~\ref{fig1}. A
  ``fuzzy'' transition temperature which would correspond more
  realistically to a mixture
  of various magnetic minerals in rock would need a general
  nonlinear function instead of $\theta\mapsto\theta\CC(\mm)$ in
  \eqref{free}. This would however make
  the specific heat $-\theta\psi_{\theta\theta}''$ dependent also on $\mm$
  and the analysis more complicated.
\end{remark}

\begin{remark}[{\sl Compressible variant with varying mass density}]\label{compres}\upshape
In situations where mass density $\varrho$ varies in space and time,
the system \eqref{MG} is to be complete by the continuity equation
$\DT\varrho=-\varrho{\rm div}\,\vv$ while the force
$-\varrho({\rm div}\,\vv)\vv/2$ in \eqref{MG-1} is to be omitted.
 The model \eqref{MG} itself used the simplification based on the assumption
 of a constant mass density $\varrho$, cf.\ the calculus \eqref{formula1},
   and neglects its variations during
   volumetric deformation (which is typically indeed small in
   liquids and in solids, too, in contrast to gases).
   Keeping the energy balance without the continuity equation
    then needs a compensation by the ``structural'' force
   $-\varrho({\rm div}\,\vv)\vv/2$, invented by R.\,Temam \cite{Tema69ASEN}
   rather for numerical purposes. This extra force is presumably very small
   (as ${\rm div}\,\vv$ is typically small) and slightly violates 
    Galilean invariancy, as pointed out in \cite{Toma21ITST}, which is
    the price for simplifying the model and its analysis. Here, in fact,
    the analysis would however work for the full model 
   because the concept of nonlinear nonsimple material ensures $\Nabla\vv$
bounded in space, although a lot of very sophisticated technicalities and
nontrivial arguments would still be needed.
\end{remark}

\begin{remark}[{\sl Convective transport of magnetization vectors}]\label{trans}
\upshape
Based on the concept of micropolar fluids \cite{Eri66TMF}, 
an isothermal model of transport of magnetization vector by the time
derivative $\DT\mm-{\bm\omega}\times\mm$ uses an angular velocity
${\bm\omega}$, cf.\
\cite{GruWei19FITM,NoSaTo16EFMN,Rose87MF,Rose02BEMF,Scro19GWPC}. It is 
possibly to be approximated as ${\rm curl}\,\vv/2$. Our choice Zaremba-Jaumann
objective transport \eqref{ZJ-vector} actually used
${\bm W}={\rm skw}(\Nabla\vv)$ instead of
$\frac12({\rm curl}\,\vv)\times\cdot$. 
A completely different approach can use rotation of magnetization vector
described in Lagrangian (reference) configuration \cite{RouTom13TMPE} or
also \cite[Sect.\,4.5.4]{KruRou19MMCM}, which would easily allow for varying
mass density but brings a lot of other technicalities, and is not well fitted
with fluids and thus with solid-fluid transition.
\end{remark}

\section{The analysis by a time-discrete approximation}\label{sec-anal}

We will consider an initial-value problem for the boundary-value problem
\eqref{MG}--\eqref{BC} considered with $\RR:=\ZJ\Ep$ as in Remark~\ref{rem-R}
by imposing the the initial conditions
\begin{align}\label{IC}
\vv|_{t=0}^{}=\vv_0,\ \ \ \ \ {\Ee}|_{t=0}^{}={\Ee}_0,\ \ \ \ \ 
  \mm|_{t=0}^{}=\mm_0,\ \ \text{ and }\ \ \W|_{t=0}^{}= \W_0\,.
\end{align}

We will use the standard notation concerning the Lebesgue and the Sobolev
spaces, namely $L^p(\varOmega;\R^n)$ for Lebesgue measurable functions
$\varOmega\to\R^n$ whose Euclidean norm is integrable with $p$-power, and
$W^{k,p}(\varOmega;\R^n)$ for functions from $L^p(\varOmega;\R^n)$ whose
all derivative up to the order $k$ have their Euclidean norm integrable with
$p$-power. We also write briefly $H^k=W^{k,2}$. The notation
$p^*$ will denote the exponent from the embedding
$W^{1,p}(\varOmega)\subset L^{p^*}(\varOmega)$, i.e.\ $p^*=dp/(d{-}p)$
for $p<d$ while $p^*\ge1$ arbitrary for $p=d$ or $p^*=+\infty$ for $p>d$.
Moreover, for a Banach space
$X$ and for $I=[0,T]$, we will use the notation $L^p(I;X)$ for the Bochner
space of Bochner measurable functions $I\to X$ whose norm is in $L^p(I)$
while $W^{1,p}(I;X)$ denotes for functions
$I\to X$ whose distributional derivative is in $L^p(I;X)$. Also, $C_{\rm w}(I;X)$
will denote the Banach space of weakly continuous functions $I\to X$, and
$C_{\rm w*}(I;X)$ of weakly* continuous if $X$ has a predual, i.e.\ there is
$X'$ such that $X=(X')^*$ where $(\cdot)^*$ denotes the dual space.
Occasionally, we will
use $L_{\rm w*}^p(I;X)$ the space of weakly* measurable mappings $I\to X$;
recall that $L_{\rm w*}^p(I;X)=L^p(I;X)$ if $X$ is separable reflexive.
Eventually, $C^1(\cdot)$ will stand for the space of
continuously differentiable functions.

The philosophy of tuning the assumptions is quite peculiar and worth
articulating the main points:
\begin{itemize}
\vspace*{-.5em}\item $p>d$ so that the velocity field gradient $\Nabla\vv$
will be surely in $L^p(I;L^\infty(\varOmega;\R^{d\times d}))$ and the transport
of $\Ee$ will be qualitatively controlled by initial conditions,
\vspace*{-.7em}\item
due to the estimation of the convective term $\vv{\cdot}\Nabla\theta$ in
\eqref{gamma-calculus} below, we need the heat capacity to be bounded, so
that the temperature gradient $\Nabla\theta$ will have only the basic
expected regularity, cf.\ \eqref{temp-w-cont} below, and
\vspace*{-.7em}\item
the exponent $\rexp$ controlling the coercivity of $\zeta(\theta;\cdot)$ and
the thus integrability of $\DT\mm$ is large enough to make the term
$\theta\CCC'(\mm)\DT\mm$ integrable when covering the standard
Landau ansatz \eqref{psi-special}.
\end{itemize}
More specifically, we will
assume, with some $C\in\R$ and some $\epsilon>0$ arbitrarily small, that
\begin{subequations}\label{ass}\begin{align}\nonumber
    &\varphi\in C^1(\R^{d\times d}{\times}\R^d),\ \ \
    \varphi(\Ee,\mm)=\FF(\Ee,\mm)+\GG(\mm)\ \ \text{ with \ $\FF\in
C^1(\R^{d\times d}{\times}\R^d)$ \ convex}, 
    \\&\nonumber\qquad\quad
    \forall(\Ee,\mm)\in\R_{\rm sym}^{d\times d}\times\R^d:
    \ \ \varphi(\Ee,\mm)\,\le\, C(1{+}|\Ee|^{2^*-\epsilon}{+}|\mm|^{2^*-\epsilon})\,,
\\&\nonumber\qquad\qquad\hspace*{9.8em}
|\varphi_\Ee'(\Ee,\mm)|\le C(1{+}|\Ee|^{2^*-1-\epsilon}{+}|\mm|^{2^*-1-\epsilon})\,,
\\&\nonumber\qquad\qquad\hspace*{9.8em}
|\varphi_\mm'(\Ee,\mm)|\le C(1{+}|\Ee|^{2^*/\rexp'}{+}|\mm|^{2^*/\rexp'})\,,
\\&\label{ass:1-}\qquad\qquad\hspace*{9.8em}
|\GG\,'(\mm)|\le C(1{+}|\mm|^{2^*/\rexp'})\,,
\\&\label{ass:1---}
\varphi\text{ is coercive:}\ \ \ \ \ \ \ \ 
\varphi(\Ee,\mm)\ge\epsilon|\Ee|^2+\epsilon|\mm|^2-C\,,
\\&\CC\in C^1(\R^d),\ \ \forall\mm,\in\R^d:\ \ \
|\CC(\mm)|\le C(1{+}|\mm|^2)\ \ \text{ and }\ \  
|\CC'(\mm)|\le C(1{+}|\mm|)\,,
\label{ass:1--}
\\&\nonumber
\phi\in C^1(\R)\,,\ \ \phi(0)=0\,,\ \ \forall\theta\,{\in}\,\R^+{:}\ \
|\phi(\theta)|\le C(1{+}\theta^{1+\epsilon})\,,
     \\&\qquad\quad\label{ass:1}
     \text{$\gamma:\theta\mapsto\phi(\theta){-}\theta\phi'(\theta)$ 
  and $\gamma^{-1}(\cdot)$ have at most linear growth,}
       \\&\label{GM-ass}
   \GM:\R\to\R\ \text{ and }\ \KK:\R\to\R
   \text{  continuous, bounded, }\inf\GM(\cdot)>0\text{ and }
\inf\KK(\cdot)>0,
   \\&\nonumber
\zeta:\R\times\R^d\to\R\
\text{ continuous,}
\\&\nonumber\qquad
\ \ \zeta(\theta;\cdot):\R^d\to\R
\ \text{ strictly convex and smooth on $\R^d{\setminus}\{0\}$,}
\\&\qquad\quad
\forall(\rr,\theta)\in\R^d\times\R^+:\ \ \
\epsilon|\rr|^\rexp\le\zeta(\theta;\rr)\le
\big(1+\big|\rr|^\rexp\big)/\epsilon\,,
\label{ass:2}
\\\label{ass:4}
&r\in\begin{cases}(30/7,6]&\text{if }d=3,
\\[-.2em]
(2,+\infty)&\text{if }d=2,
\end{cases}
\ \ p>\rexp,\ \
\varrho,\,\kappa,\,\varkappa,\,
\DD,\,\NU,\,\mu_0>0\,,
\\
&\vv_0\!\in\! L^2(\varOmega;\R^d),\ \ \
        {\Ee}_0\!\in\! H^1(\varOmega;\R_{\rm sym}^{d\times d}),\ \ \
 \label{ass:5}
 \mm_0\!\in\! H^1(\varOmega;\R^d),\ \ \
 \theta_0\!\in\! L^1(\varOmega),\ \ \ \theta_0\ge0,\ \,
\\\label{ass:6}&\ff\in L^1(I;L^2(\varOmega;\R^d))\,,\ \ \ \hh_\text{\rm geo}\in
L^1(I;H^1(\varOmega;\R^d))\,,
\ \ \ j_{\rm ext}^{}\in L^1(I{\times}\varGamma)\,,\ \ j_{\rm ext}^{}\ge0\,.
\end{align}\end{subequations}
Let us note that we do not assume $\GG$ and $\CC$ convex in order to cover
the Landau magnetic phase-transition model \eqref{psi-special},
using particularly $\GG(\mm)=-a_0\theta_{\rm c}|\mm|^2$,
and the split $\varphi(\Ee,\mm)=\FF(\Ee,\mm)+\GG(\mm)$ is related
with the time-discretization used in the proof rather than with
the problem itself. For this, we will need to make a
regularization of the magnetization flow-rule when discretized in time,
cf.\ the magnetic flow rule \eqref{ED-4+disc} and then also the
(\ref{ED+disc}a,f) below.
Note also that \eqref{ass:1---} admits a canonical
choice $\phi(\theta)=c_{\rm v}\theta(1{-}{\rm ln}\theta)$, which gives
$\gamma(\theta)=c_{\rm v}\theta$ with a constant heat capacity $c_{\rm v}$.

Since $\Nabla^2u$ has a sense as an distribution only, for the weak formulation
\eqref{MG-1-weak} of the momentum equation \eqref{MG-1} below, we rather
use the Green formula:
\begin{align}
\int_\varOmega\widetilde\vv{\cdot}(\Nabla^2u)\mm\,\d x
=\int_\varGamma(\!\!\!\!\lineunder{\widetilde\vv{\cdot}\nn}{$=0$}\!\!\!\!)(\Nabla u{\cdot}\mm)\,\d x
-\int_\varOmega(\Nabla u{\cdot}\mm){\rm div}\,\widetilde\vv
+(\widetilde\vv{\otimes}\Nabla u){:}\Nabla\mm\,\d x\,.
\label{Green-for-momentum}\end{align}

\begin{definition}[Weak solutions]\label{def}
 A six-tuple $(\vv,\Ee,\RR,\mm,u,\theta)$ with
\begin{subequations}\begin{align}
 &\vv\!\in\! C_{\rm w}(I;L^2(\varOmega;\R^d))\cap L^p(I;W^{2,p}(\varOmega;\R^d))
    \nonumber\\[-.2em]&\hspace*{9.5em}
    \ \mbox{ with \ $\frac{\partial}{\partial t}$}
    \vv\!\in\!L^{p'}\!(I;W^{2,p}(\varOmega;\R^d)^*)+L^1(I;L^2(\varOmega;\R^d))\,,
\\&\Ee\in C_{\rm w}(I;H^1(\varOmega;\R_{\rm sym}^{d\times d}))\,\cap\,H^1(I,L^2(\varOmega;\R_{\rm sym}^{d\times d}))\,,
    \\& \RR\in L^2(I;H^1(\varOmega;\R_{\rm dev}^{d\times d}))\,,
 \\&\mm\in C_{\rm w}(I;H^1(\varOmega;\R^d))\ \mbox{ with \ }
 \ZJ\mm\in L^\rexp(I{\times}\varOmega;\R^d)\ \text{ and }\ \Delta\mm
\in L^{\rexp'}(I{\times}\varOmega;\R^d)\,,
\\&u\in C_{\rm w}(I;H^1(\R^d))
\;\ \ \ \ \ \text{ with }\ \
\mbox{$\frac{\partial}{\partial t}$}u\big|_{I{\times}\varOmega}^{}\in L^2(I;H^1(\varOmega))\,,
\\&\theta
    \in C_{\rm w}(I;L^1(\varOmega))\,\cap\,L^\sigma(I;W^{1,\sigma}(\varOmega))\ 
    \mbox{ with $1\le\sigma<\frac{d{+}2}{d{+}1}$},
      \ \ \text{and}\ \theta\ge0\ \text{ a.e. on }I{\times}\varOmega\,,
    \label{temp-w-cont}
 \end{align}\end{subequations}
 will be called a weak solution to the boundary-value problem
 \eqref{MG}--\eqref{BC} with the initial conditions \eqref{IC} if
${\bm S}_\text{\sc e}=\varphi_\Ee'(\Ee,\mm){+}\psi(\Ee,\mm,\theta)\bbI
\in L^\infty(I;L^1(\varOmega;\R_{\rm sym}^{d\times d}))$,
 $\varphi_\mm'(\Ee,\mm)+\theta\CC'(\mm)
\in L^2(I{\times}\varOmega;\R^d)$,
   $\nn{\cdot}\vv=0$ on $I{\times}\varGamma$,
   $\W=\gamma(\theta)\in C_{\rm w}(I;L^1(\varOmega))$
    with $\frac{\partial}{\partial t}\W\in
    L^1(I;H^{d+1}(\varOmega)^*)$, and if

\vspace*{-2.2em}

  \begin{subequations}\label{def+}\begin{align}\nonumber
      &\int_0^T\!\!\!\int_\varOmega\bigg(\varrho\Big(\!(\vv{\cdot}\Nabla)\vv
      {+}\frac12({\rm div}\,\vv)\vv\Big){\cdot}\widetilde\vv
      +\big({\bm S}_\text{\sc e}
      {+}\DD\EE(\vv)\big){:}\EE(\widetilde\vv)
          +{\bm S}_\text{\sc c}{:}\Nabla\widetilde\vv 
        \\[-.2em]&\hspace{4em}\nonumber
+\NU|\Nabla\EE(\vv)|^{p-2}\Nabla\EE(\vv)\Vdots\Nabla\EE(\widetilde\vv)
+\mu_0(\Nabla\hh_\text{\rm geo})^\top\mm
{\cdot}\widetilde\vv
-\mu_0(\Nabla u{\cdot}\mm){\rm div}\,\widetilde\vv
 \\[-.2em]&\hspace{4em}
-\mu_0(\widetilde\vv{\otimes}\Nabla u){:}\Nabla\mm
-\varrho\vv{\cdot}\pdt{\widetilde\vv}\bigg)\,\d x\d t
=\int_\varOmega\varrho\vv_0{\cdot}\widetilde\vv(0)\,\d x
      +\int_0^T\!\!\!\int_\varOmega\ff{\cdot}\widetilde\vv\,\d x\d t
 \label{MG-1-weak}
 \intertext{with ${\bm S}_\text{\sc c}\in L^1(I{\times}\varOmega;\R^{d\times d})$
        from \eqref{MG-1}
        for all $\widetilde\vv\in C^1(I{\times}\varOmega;\R^d)
        \cap L^p(I;W^{2,p}(\varOmega;\R^d))$
        with $\nn{\cdot}\widetilde\vv=0$ on $I{\times}\varGamma$ and
        $\widetilde\vv(T)=0$, and if}
&\nonumber\\[-2.3em]
&\int_0^T\!\!\!\int_\varOmega
\big(\GM(\W)\RR{+}{\rm dev}\,\bm{S}\big){:}\widetildeEp
+\varkappa\Nabla\RR\Vdots\Nabla\widetildeEp
\,\d x\d t=0
\label{MG-2-weak}
\intertext{holds for all $\widetildeEp\in
  H^1(I{\times}\varOmega;\R_{\rm dev}^{d\times d})$,
  and if}\nonumber
     &\int_0^T\!\!\!\int_\varOmega\bigg(\zeta(\theta;\widetilde\rr)
 +\Big(\hh_\text{\rm geo}-\frac{\varphi_\mm'(\Ee,\mm)+
 \theta\CC'(\mm)}{\mu_0}
-\nabla u\Big){\cdot}\big(\widetilde\rr{-}\ZJ\mm\big)
\\[-.6em]&\hspace*{11.5em}
+\kappa\Nabla\mm{:}\Nabla\widetilde\rr+\kappa\Delta\mm
{\cdot}\big((\vv{\cdot}\Nabla)\mm
-{\rm skw}(\Nabla\vv)\mm\big)\bigg)\,\d x\d t
\nonumber\\[-.5em]&\hspace*{8em}
      +\int_\varOmega\frac{\kappa}2|\Nabla\mm_0|^2\,\d x\ge\int_\varOmega
    \frac{\kappa}2|\Nabla\mm(T)|^2\,\d x
    +\int_0^T\!\!\!\int_\varOmega\zeta\big(\theta;\ZJ\mm\big)\,\d x\d t
   \label{VI}
\intertext{holds for all $\widetilde\rr\in L^2(I;H^1(\varOmega;\R^d))$, and
at all time instants}
&\int_{\R^d}\nabla u{\cdot}\nabla\widetilde u\,\d x
=\int_\varOmega\mm{\cdot}\nabla\widetilde u\,\d x
\label{MG-4-weak}
\intertext{holds for all $\widetilde u\in H^1(\R^d)$, and}
&\nonumber
\int_0^T\!\!\!\int_\varOmega\!
\bigg(\big(\KK(\theta)\nabla\theta{-}\W\vv\big){\cdot}
\nabla\widetilde\W-\Big(\xi\big(\theta;\EE(\vv),\RR,\ZJ{\mm}\big)
+\theta\CCC'(\mm){\cdot}\DT\mm+
\big(\theta\CCC(\mm){+}\phi(\theta)\big)\,{\rm div}\,\vv\Big)\widetilde\W
\\[-.4em]&\label{ED-6+weak}
\hspace*{12em}
-\W\pdt{\widetilde\W}\bigg)
\,\d x\d t
=\int_\varOmega\gamma(\theta_0)\widetilde\W(0)\,\d x
+\int_0^T\!\!\!\int_\varGamma j_{\rm ext}\widetilde\W\,\d S\d t
\end{align}\end{subequations} 
with $\xi(\theta;\EE(\vv),\RR,\ZJ{\mm})$ from \eqref{MG-6+thermo}
for any $\widetilde\W\in W^{1,\infty}(I{\times}\varOmega)$ with
$\widetilde\W(T)=0$, and eventually also
$\ZJ\Ee=\EE(\vv)-\RR$ holds a.e.\ on $I{\times}\varOmega$
and $\Ep(0)=\Ep_0$ and  $\Ee(0)=\Ee_0$ a.e.\ on $\varOmega$.
\end{definition}

\begin{proposition}[Existence of weak solutions]
Let \eqref{ass} be valid. Then the initial-boundary-value problem for the
system \eqref{MG} with the boundary and initial conditions \eqref{BC} and
\eqref{IC} has a weak solution $(\vv,\Ee,\RR,\mm,u,\theta)$ according
Definition~\ref{def} and every such solution also satisfies
the mechanical-energy balance \eqref{energy+} and
conserves the total energy in the sense \eqref{energy-mag}.
\end{proposition}

\def\Seetauk{{\bm S}_{\text{\sc e},\etau}^k}
\def\Svetauk{{\bm S}_{\text{\sc v},\etau}^k}

\noindent{\it Proof.}  For lucidity, we divide the proof into six steps.

\medskip\noindent{\it Step 1: Approximation by time discretization.}
As we need testing by convective (but not mere partial) time derivatives, using
of Galerkin method would be very technical (if not just impossible). Therefore, 
we use the Rothe method, i.e.\ the fully implicit time discretization with an
equidistant partition of the time interval $I$ with the time step $\tau>0$.
This approximation is also rather technical because
the stored energy $\varphi(\Ee,\cdot)+\theta\CC(\cdot)$
is intentionally nonconvex for $\theta<\theta_{\rm c}$,
and its test by the time difference of $\mm$ is to be estimated ``on the
right-hand side''. For this reason, we denote the possibly nonconvex part
of the stored energy by
$$
\varpi(\mm,\theta):=\GG(\mm)+\theta\CC(\mm)
$$
for this proof. In addition to time discretization, we thus use also
a modification of $\varpi$ by some
$\varpi_\eps(\mm,\theta)=\GG_\eps(\mm)+\theta\CC_\eps(\mm)$
 with $\GG_\eps,\CC_\eps\in W^{1,\infty}(\R^d)$, specifically we put
\begin{align}\label{regul}
\GG_\eps(\mm):=\frac{\GG(\mm)}{1{+}\eps|\mm|^2}
\ \ \text{ and }\ \ 
\CC_\eps(\mm):=\frac{\CC(\mm)}{1{+}\eps|\mm|^2}\,,
\ \ \text{ so that }\ \ 
[\varpi_\eps]_\theta'(\mm,\theta):=\CC_\eps(\mm)\,.
\end{align}
Consequently, we also put $\psi_\eps(\Ee,\mm,\theta)=\FF(\Ee,\mm)
+\GG_\eps(\mm)+\theta\CC_\eps(\mm)+\mu_0\hh_{\rm dem}{\cdot}\mm-\phi(\theta)$.

We denote by $\vvk$, $\Eek$, ${\bm S}_\etau^k$, ... the approximate
  values of $\vv$, $\Ee$, ${\bm S}$, ...
  at time $k\tau$ with $k=1,2,...,T/\tau$.
  We introduce a shorthand notation for the bi-linear operators
 \begin{subequations}\begin{align}
 &\bm b_\text{\sc zj}^{}(\vv,\mm)=(\vv{\cdot}\Nabla)\mm
  -{\rm skw}(\Nabla\vv)\mm\quad\text{ and}
      \\&
      \bm B_\text{\sc zj}^{}(\vv,\bm E)=(\vv{\cdot}\Nabla)\bm E
      -{\rm skw}(\Nabla\vv)\Ee+\Ee\,{\rm skw}(\Nabla\vv) \,.
    \end{align}
\end{subequations}
We devise a fully coupled time discretization; here let us remark
that usual efficient decoupled (staggered) schemes would be problematic
in this case because the convective corrotation terms and the adiabatic
terms bond intimately the particular equations and their decoupling would
complicate the a-priori estimates. On the other hand, we easily can
use delayed temperature in some dissipation terms, cf.\ (\ref{ED+disc}c,d,f)
below. Beside the mentioned regularization of $\CC$, we use a bit lower
dissipation heat source in \eqref{ED-6+disc} by using the factor $1-\eps$.
We will then use the following recursive regularized time-discrete scheme
written in the classical form as 
\begin{subequations}\label{ED+disc}\begin{align}
      \nonumber
  &\varrho \Big(\frac{\vvk{-}\vvkk\!\!\!}\tau+
 (\vvk{\cdot}\Nabla)\vvk\Big)
 ={\rm div}\big(\Seetauk+\Svetauk+{\bm S}_{\text{\sc c},\etau}^k\!
-\mu_0((\hh_\etau^k){\cdot}\mmk)\bbI-{\rm div}\mathfrak{H}_\etau^k\big)
\\\nonumber
  &\hspace*{18em}
-\frac\varrho2({\rm div}\,\vvk)\,\vvk
+\mu_0(\Nabla\hh_\tau^k)^\top\mmk
+\varrho{\bm g}_\tau^k\,
 \\\nonumber
  &\hspace*{1em}\text{ with }\ \ 
     \Seetauk=\varphi_{\Ee}'(\Eek,\mmk){+}\psi_\eps(\Eek,\mmk,\theta_\etau^k)\bbI\,,
     \ \ \hh_\etau^k=\hh_{\text{\rm geo},\etau}^k{-}\nabla u_\etau^k\,,
 \\\nonumber
  &\hspace*{4.5em}\Svetauk=\DD\EE(\vvk)
     \,,\ \ \ \mathfrak{H}_\etau^k=\NU|\Nabla\EE(\vvk)|^{p-2}\Nabla\EE(\vvk)\,,
\ \ \text{ and }\ \ \
\\&\hspace*{4.5em}
  {\bm S}_{\text{\sc c},\etau}^k\!=\kappa\mu_0
  \Nabla\mmk{\otimes}\Nabla\mmk
  -\frac{\kappa\mu_0}{\qexp}|\Nabla\mmk|^{\qexp}\bbI
-\mu_0{\rm skw}\big(
   (\hh_{\rm drv,\etau}^k){\otimes}\mmk\big)\,,\label{ED-1+disc}
\\[-.3em]&\frac{\bm E_\etau^k{-}\bm E_\etau^{k-1}\!\!\!}\tau
+\bm B_\text{\sc zj}^{}(\vvk,\bm E_\etau^k)=\EE(\vvk)-\RRk\,,
\label{ED-2+disc}
 \\\label{ED-3+disc}
&
\GM(\theta_\etau^{k-1})\RR_\etau^k+{\rm dev}\,\Seetauk=\varkappa\Delta\RR_\etau^k\,,
\\\nonumber
&
\partial_{\rr}\zeta\big(\theta_\etau^{k-1};\rrk\big)\ni
\hh_{\rm drv,\etau}^k\ \ \ \text{ with }\ \ \
\rrk=\frac{\mmk{-}\mm_\etau^{k-1}\!\!\!}\tau
+\bm b_\text{\sc zj}^{}(\vvk,\mm_\etau^k)\ \ \text{and }\ \ 
\\[-.1em]\label{ED-4+disc}
&\hspace*{5em}
\hh_{\rm drv,\etau}^k=\hh_{\text{\rm geo},\tau}^k
-\frac{\FF_\mm'(\Ee_\etau^k,\mmk){+}[\varpi_\eps]_\mm'(\mmk,\thetak)\!}{\mu_0}
-\Nabla u_\etau^k+\kappa\Delta\mmk\,,
\\&
\Delta u_\etau^k={\rm div}(\chi_\varOmega^{}\mm_\etau^k)\ \ \ \
\text{ on }\ \ \R^d\,,
\label{ED-5+disc}\\
\nonumber
&\frac{\W_\etau^k{-}\W_\etau^{k-1}\!\!\!}\tau
+
{\rm div}\big(\W_\etau^k\vvk-\KK(\theta_\etau^{k-1})\Nabla\thetak\big)=
(1{-}\eps)\Big(\GM(\theta_\etau^{k-1})|\RR_\etau^k|^2+\DD|\EE(\vvk)|^2
\\[-.1em]\nonumber
&\hspace*{4em}
+\NU|\Nabla\EE(\vvk)|^p+\partial_{\rr}\zeta\big(\theta_\etau^{k-1};
\rrk\big){\cdot}\rrk+\varkappa|\Nabla\RR_\etau^k|^2\Big)
\\[-.0em]\nonumber
&\hspace*{4em}
+\thetak\CCC_\eps'(\mmk){\cdot}\Big(\frac{\mmk{-}\mmkk\!\!\!}\tau
 +(\vvk{\cdot}\Nabla)\mmk\Big)
\\&\hspace*{4em}
+\big(\thetak\CCC_\eps(\mmk){+}\phi(\thetak)\big){\rm div}\,\vvk
\hspace*{3em}\text{ with }\ \ \W_\etau^k=\gamma(\thetak)\,.
 \label{ED-6+disc}
 \end{align}\end{subequations}
We complete the system (\ref{ED+disc}a-d,f) by the corresponding boundary
conditions, i.e.
\begin{subequations}\label{BC-disc}
  \begin{align}\label{BC1-disc}
&\!\big[\big(\Seetauk+\Svetauk+{\bm S}_{\text{\sc c},\etau}^k\!
               -\mu_0(\hh_\etau^k{\cdot}\mmk)\bbI
  -{\rm div}\mathfrak{H}_\etau^k\big)\nn+\divS(\mathfrak{H}_\etau^k\nn)
  \big]_\text{\sc t}^{}=0\,,
  \\\label{BC2-disc}
   &\vvk{\cdot}\nn=0\,,\ \ \ \ \Nabla\EE(\vvk){:}(\nn{\otimes}\nn)=0\,,\ \ \ 
  (\nn{\cdot}\Nabla)\RR_\etau^k=0\,,\ \ \ 
 (\nn{\cdot}\Nabla)\mmk=0,\ \ \text{ and }\ \
 \\[-.0em]
  &\nn{\cdot}\KK(\theta_\etau^{k-1})\Nabla\thetak=j_{\rm ext,\tau}^k
  \,,\label{BC3-disc}
\end{align}
\end{subequations}
while \eqref{ED-5+disc} completes with the condition
$\lim_{|x|\to\infty}u_\etau^k(x)=0$ for all $k=1,...,T/\tau$.
Here we used $\ff_\tau^k:=\int_{(k-1)\tau}^{k\tau}\ff(t)\,\d t$ and
similarly also for $\hh_{\text{\rm geo},\tau}^k$ and $j_{\rm ext,\tau}^k$.
This system of boundary-value problems is to be solved
recursively for $k=1,2,...,T/\tau$, starting 
with the initial conditions for $k=1$: 
\begin{align}\label{IC-disc}
\vv_\etau^0=\vv_0,\ \ \ \ \ {\Ee}_\etau^0={\Ee}_0,\ \ \ \ \ 
  \mm_\etau^0=\mm_0,\ \ \text{ and }\ \ \W_\etau^0= \W_0 \,.
\end{align}

The existence of a weak solution $(\vvk,\Eek,
\RR_\etau^k,\mmk,u_\etau^k,\W_\etau^k)\in W^{2,p}(\varOmega;\R^d)\times
L^2(\varOmega;\R_{\rm sym}^{d\times d})
\times H^1(\varOmega;\R_{\rm dev}^{d\times d})\times H^1(\varOmega;\R^d)\times H^1(\R^d)\times W^{1,1}(\varOmega)$ 
of the coupled quasi-linear boundary-value problem
(\ref{ED+disc})--(\ref{BC-disc}) can thus be seen by a combination
of the quasilinear technique for \eqref{ED-1+disc} involving the quasilinear
term ${\rm div}^2(\NU|\Nabla e(\vvk)|^{p-2}\Nabla e(\vvk))$, with the usual
semi-linear technique for (\ref{ED+disc}c--d), with the $L^1$-technique for the
heat equation \eqref{ED-6+disc}, and with the set-valued inclusion (in fact a
variational inequality) \eqref{ED-4+disc} provided $\vvkk$, $\Eekk$, 
$\mmkk$, and $\W_\etau^{k-1}$ are known from the previous time step. Actually,
$\W_\etau^k\in W^{1,q}(\varOmega)$ for any $q<d'$. Let us note that the system
\eqref{ED+disc} is indeed fully coupled due to the convective derivatives
and due to the adiabatic effects, cf.\ the term
$\phi(\thetak)$ in \eqref{ED-1+disc}, and it seems
not possible to devise some decoupled (staggered) discrete scheme which
would allow for some reasonable estimation strategy of the recursive scheme.
The mentioned coercivity is a particular consequence of the a-priori estimates
derived below. Thus also $\RR_\etau^k\in H^1(\varOmega;\R_{\rm dev}^{d\times d})$ is
obtained. Let us note that, due to the convective terms,
this system does not have any potential so the rather nonconstructive Brouwer
fixed-point arguments combined with the Galerkin approximation. Also strict
monotonicity of the main parts of  (\ref{ED+disc}a--d) so that the approximated
right-hand side of the semilinear equation \eqref{ED-6+disc} converges strongly
in $L^1(\varOmega)$. In general, one cannot expect any uniqueness of this
solution.

Moreover, $\W_\etau^k\ge0$ a.e.\ on $\varOmega$ for at least one weak solution.
To show it, we test \eqref{ED-6+disc} by the negative part of $\thetak$.
Here we exploit $\phi(0)=0$ so that the adiabatic terms vanish
for $\theta\to0+$, and thus prove that $\thetak\ge0$. Using also
$\gamma([0,+\infty))\ge0$, then also $\W_\etau^k=\gamma(\thetak)\ge0$.

Using the values $(\vvk)_{k=0}^{T/\tau}$, we define the piecewise constant and
the piecewise affine interpolants respectively as
\begin{align}\label{def-of-interpolants}
&\overlinevvtau(t):=\vvk,\ \ \ \underline\vv_\etau(t):=\vvkk,
\ \text{ and }\ \vv_\etau(t):=\Big(\frac t\tau{-}k{+}1\Big)\vvk
\!+\Big(k{-}\frac t\tau\Big)\vvkk
\end{align}
for $(k{-}1)\tau<t\le k\tau$ with $k=0,1,...,T/\tau$. Analogously, we define
also $\Eetau$, $\overlineEetau$, $\overlineRRtau$,  $\underline{\W}_\etau$,
etc. Thus, \eqref{ED+disc} holding a.e.\ on $\varOmega$ for $k=1,...,T/\tau$
can be written ``compactly'' as 
 \begin{subequations}\label{ED+d}\begin{align}
 \nonumber
 &\varrho
 \Big(\frac{\partial\vv_\etau}{\partial t}
 +(\overlinevvtau{\cdot}\Nabla)\overlinevvtau\Big)
 ={\rm div}\big(\overlineSetau+\overlineSvtau
     +\overlineSstrtau\!-\mu_0(\overline{\hh}_\etau{\cdot}\overlinemmtau)\bbI
               -{\rm div}\overline{\mathfrak{H}}_\etau\big)
\\\nonumber
  &\hspace*{15em}
-\frac\varrho2({\rm div}\,\overlinevvtau)\,\overlinevvtau
+\mu_0(\overline{\hh}_\etau)^\top\overlinemmtau+\varrho\overline\ff_\tau\,
 \\\nonumber
  &\hspace*{1em}\text{ with }\ \ 
     \overlineSetau=\varphi_\Ee'(\overlineEetau,\overlinemmtau){+}
     \psi_\eps(\overlineEetau,\overlinemmtau,\overline\theta_\etau)\bbI\,,\ \ \
\overline{\hh}_\etau=\overline\hh_{\text{\rm geo},\tau}{-}\nabla\overline u_\etau\,,
\\&\hspace*{4.5em}\nonumber
\overlineSvtau=\DD\EE(\overlinevvtau)
     \,,\ \ \ \overline{\mathfrak{H}}_\etau=\NU|\Nabla e(\overlinevvtau)|^{p-2}\Nabla e(\overlinevvtau)\,,\ \ \text{ and }
\\&\hspace*{4.5em}
  \overlineSstrtau\!=\kappa\mu_0
  \Nabla\overlinemmtau{\otimes}\Nabla\overlinemmtau
  -\frac{\kappa\mu_0}{\qexp}|\Nabla\overlinemmtau|^{\qexp}\bbI
   -\mu_0{\rm skw}\big(\overline\hh_{\rm drv,\etau}{\otimes}\overlinemmtau\big)\,,
\label{ED-1+d}
\\[-.3em]&
\frac{\partial\Ee_\etau}{\partial t}
+\bm B_\text{\sc zj}^{}(\overlinevvtau,\overlineEetau)
=\EE(\overlinevvtau)-\overlineRRtau\,,
 \label{ED-2+d}
 \\\label{ED-3+d}
 &
\GM(\underline\W_\etau)\overlineRRtau
+{\rm dev}\,\overlineSetau =
\varkappa\Delta\overlineRRtau\,,
\\\nonumber
&
\partial_{\rr}\zeta\big(\underline\theta_\etau;\overline\rr_\etau\big)\ni
\overline\hh_{\rm drv,\etau}\ \ \
\text{ with }\ \ \ \overline\rr_\etau=
\frac{\partial\mm_\etau}{\partial t}+
 \bm b_\text{\sc zj}^{}(\overlinevvtau,\overlinemmtau)\ \ \ \ \text{ and }
\\[-.1em]\label{ED-4+d}
&\hspace*{3.5em}
\overline\hh_{\rm drv,\etau}=\overline\hh_{\text{\rm geo},\tau}
-\frac{\FF_\mm'(\overlineEetau,\overlinemmtau){+}
[\varpi_\eps]_\mm'(\overlinemmtau,\overline\theta_\etau)
\!}{\mu_0}-\nabla\overline u_\etau+\kappa\Delta\overlinemmtau\,,
\\
\nonumber
&\frac{\partial\W_\etau}{\partial t}
+{\rm div}\big(\overline\W_\etau\overlinevvtau
-\KK(\underline\W_\etau)\Nabla\overline\theta_\etau\big)=
(1{-}\eps)\Big(\GM(\underline\theta_\etau)|\overlineRRtau|^2
+\DD|\EE(\overlinevvtau)|^2
\\[-.3em]\nonumber
&\hspace*{2em}+\NU|\Nabla\EE(\overlinevvtau)|^2+
\partial_{\rr}\zeta\big(\underline\theta_\etau;\overline\rr_\etau\big){\cdot}\overline\rr_\etau
+\varkappa|\Nabla\overlineRRtau|^2\Big)
+\overline\theta_\etau\CCC_\eps'(\overlinemmtau)
{\cdot}\Big(
\frac{\partial\mm_\etau\!\!}{\partial t}
+\overlinevvtau{\cdot}\Nabla\overlinemmtau\Big)
\\[-.2em]
&\hspace*{6em}
+\big(\overline\theta_\etau\CCC_\eps(\overlinemmtau)
{+}\phi(\overline\theta_\etau)\big){\rm div}\,\overlinevvtau
 \qquad\text{ with }\quad
\overline\W_\etau=\gamma(\overline\theta_\etau)
 \label{ED-6+d}
\intertext{
 holding on $I{\times}\varOmega$ either a.e.\ or in a weak sense
 involving also the boundary conditions \eqref{BC-disc} which are
 to be written analogously in terms of the above introduced interpolants, and}
&\label{ED-5+d}
\Delta\overline u_\etau={\rm div}(\chi_\varOmega^{}\overlinemmtau)
\end{align}\end{subequations}
holding on $I{\times}\R^d$ in the weak sense together with the condition
$\lim_{|x|\to\infty}\overline u_\etau(t,x)=0$ for a.a.\ $t\in I$.

\medskip\noindent{\it Step 2: A-priori estimates}.
The a-priori estimation is based on the energy test for the mechanical
part combined with the heat problem. This means here 
the test of \eqref{ED-1+disc} by $\vvk$ while using also \eqref{ED-2+disc}
tested by $\bm{S}_\etau^k=\varphi_\Ee'(\Eek,\mmk)$, then the test
\eqref{ED-3+disc} by $\RR_\etau^k$, the inclusion \eqref{ED-4+disc} by
$\mu_0\rrk$, \eqref{ED-5+disc} by $(u_\etau^k{-}u_\etau^{k-1})/\tau$,
and  \eqref{ED-6+disc} by 1. 
We thus obtain an energy-like inequality for the time-discrete approximation
corresponding to  \eqref{energy-mag}, but in contrast to \eqref{energy-mag},
we treat the nonconvex terms $\CC_\eps$ and $\GG_\eps$ 
``on the right-hand side''.

More specifically, the terms related to inertia in \eqref{ED-1+disc} uses
the calculus
\begin{align}\nonumber
  \Big(\varrho\frac{\vvk{-}\vvkk\!\!}\tau
 & +\varrho(\vvk{\cdot}\Nabla)\vvk
 +\frac\varrho2(\operatorname{div}\vvk)\vvk\Big){\cdot}\vvk=
\frac{\varrho}2\frac{|\vvk|^2-|\vvkk|^2\!\!}\tau
  \\[-.4em]&
   +\varrho(\vvk{\cdot}\Nabla)\vvk\cdot\vvk
 +\frac{\varrho}2({\rm div}\,\vvk)|\vvk|^2
  +\tau\frac{\varrho}2\Big|\frac{\vvk{-}\vvkk\!\!}\tau\ \Big|^2\,.
\label{test-of-convective}\end{align}
This holds pointwise and, when integrated over $\varOmega$, we further use
also \begin{align}\nonumber
 \int_\varOmega\varrho(\vvk{\cdot}\Nabla)\vvk{\cdot}\vvk\,\d x&=
  -\int_\varOmega\frac{\varrho}2|\vvk|^2({\rm div}\,\vvk)\,\d x
  +\int_\varGamma\frac{\varrho}2|\vvk|^2(\vvk{\cdot}\nn)\,\d S
\\&=
-\int_\varOmega\Big(\frac{\varrho}2|\vvk|^2\bbI\Big):\EE(\vvk)\,\d x
  +\int_\varGamma\frac{\varrho}2|\vvk|^2(\vvk{\cdot}\nn)\,\d S\,.
\label{convective-tested+}\end{align}
The last term in \eqref{test-of-convective} is non-negative and can simply be
forgotten, which will give the inequality
\begin{align}\nonumber
  &\int_\varOmega\!\Big(\varrho\frac{\vvk{-}\vvkk\!\!}\tau
  +\varrho(\vvk{\cdot}\Nabla)\vvk
  +\frac\varrho2({\rm div}\,\vvk)\vvk\Big)\cdot\vvk\,\d x
\\&\qquad\nonumber
\stackrel{\eqref{test-of-convective}}{\ge}\int_\varOmega\!\Big(
\frac{\varrho}2\frac{|\vvk|^2-|\vvkk|^2\!\!}\tau
 +\varrho(\vvk{\cdot}\Nabla)\vvk{\cdot}\vvk
  +\frac\varrho2({\rm div}\,\vvk)|\vvk|^2\Big)\,\d x
   \\&\qquad\qquad\stackrel{\eqref{convective-tested+}}{=}
  \int_\varOmega\frac{\varrho}2\frac{|\vvk|^2\!-|\vvkk|^2\!\!}\tau\,\d x
  +\int_\varGamma\frac\varrho2|\vvk|^2(\vvk{\cdot}\nn)\,\d S\,.
\label{test-of-convective+}\end{align}
The last term vanishes due to the boundary condition \eqref{BC1-disc}.
The further term in \eqref{ED-1+disc} uses the calculus
\begin{align}\nonumber
&\!\!\!\int_\varOmega(\operatorname{div}\Seetauk){\cdot}\vvk\,\d x=\int_\varGamma\Seetauk
  {:}(\vvk{\otimes}\nn)\,\d S-\int_\varOmega\Seetauk{:}\EE(\vvk)\,\d x
\\&\nonumber\ \ \stackrel{\eqref{ED-2+disc}}{=}
\int_\varGamma\Seetauk{:}(\vvk{\otimes}\nn)\,\d S
-\int_\varOmega\Seetauk{:}\Big(\frac{\Ee_\etau^k{-}\Ee_\etau^{k-1}\!\!}\tau+
\bm B_\text{\sc zj}^{}(\vvk,\Ee_\etau^k)-\RR_\etau^k\Big)\,\d x
\\&\nonumber\ \stackrel{(\ref{ED+disc}a,c)}{=}\!\!\int_\varGamma\Seetauk
{:}(\vvk{\otimes}\nn)\,\d S
+\int_\varOmega\bigg(
\GM(\theta_\etau^{k-1})|\RR_\etau^k|^2+\varkappa|\Nabla\RR_\etau^k|^2
-\big(\psi_\eps(\Eek,\mmk,\theta_\etau^k)\bbI
\\[-.4em]&\hspace*{13em}+\varphi_\Ee'(\Eek,\mmk)\big)
    {:}\Big(\frac{\Eek{-}\Eekk\!\!}\tau
    +\bm B_\text{\sc zj}^{}(\vvk,\Eek) \Big)\bigg)\,\d x\,,
\label{test-by-v}\end{align}
where we used also \eqref{ED-3+disc} tested by $\RR_\etau^k$.

Further, \eqref{ED-4+disc} is to be tested by $\mu_0\rrk$ and use an
analog of the calculus \eqref{test-damage}--\eqref{formula3} by 
exploiting the algebra
\begin{align}
\hh_{\text{\rm geo},\tau}^k{\cdot}\frac{\mmk{-}\mmkk\!\!}\tau
=\frac{\hh_{\text{\rm geo},\tau}^k{\cdot}\mmk-\hh_{\text{\rm geo},\tau}^{k-1}{\cdot}\mmkk\!\!}\tau
-\frac{\hh_{\text{\rm geo},\tau}^k{-}\hh_{\text{\rm geo},\tau}^{k-1}\!\!}\tau{\cdot}\mmkk
\end{align}
and convexity of the functional
$\mm\mapsto\int_{\R^d}\frac1\qexp\kappa|\Nabla\mm|^\qexp\,\d x$. Thus we will get
the inequality as a discrete analog of \eqref{formula3}, i.e.
\begin{align}\nonumber
&\int_\varOmega\mu_0\hh_{\rm drv,\etau}^k{\cdot}\rrk\,\d x
\le\int_\varOmega\!
\mu_0\frac{\hh_{\text{\rm geo},\tau}^k{\cdot}\mmk-\hh_{\text{\rm geo},\tau}^{k-1}{\cdot}\mmkk\!\!}\tau
-\mu_0\kappa\frac{|\Nabla\mmk|^\qexp-|\Nabla\mmkk|^\qexp\!\!}{\qexp\tau}\,\d x
\\[-.2em]&\nonumber\
+\!\int_\varOmega\bigg(\frac{\mu_0\kappa}\qexp|\Nabla\mmk|^\qexp{\rm div}\,\vvk
-\mu_0\kappa(\Nabla\mmk{\otimes}\Nabla\mmk){:}\EE(\vvk)
-\FF_{\mm}'(\Eek,\mmk){\cdot}\frac{\mmk{-}\mmkk\!\!}{\tau}
\\[-.3em]&\nonumber\
+\big(\hh_{\text{\rm dem},\tau}^k{-}[\varpi_\eps]_\mm'(\mmk,\thetak)\big){\cdot}
\Big(\frac{\mmk{-}\mmkk\!\!}{\tau}+(\vvk{{\cdot}}\Nabla)\mmk\Big)
-\mu_0\frac{\hh_{\text{\rm geo},\tau}^k{-}\hh_{\text{\rm geo},\tau}^{k-1}\!\!}\tau{\cdot}\mmkk
\\[-.3em]&\
-\mu_0(\hh_{\text{\rm geo},\tau}^k{\cdot}\mmk){\rm div}\,\vvk\!
-\mu_0(\nabla\hh_{\text{\rm geo},\tau}^k)^{\top}\mmk{\cdot}\vvk\!
-\mu_0{\rm skw}\big(\hh_{\rm drv,\etau}^k{\otimes}\mmk\big)
{:}\Nabla\vvk\bigg)\,\d x\,.
\end{align}
For this test, we also exploit that
$\partial_{\rr}\zeta(\theta;\rr){\cdot}\rr$ is uniquely defined
even though the subdifferential $\partial_{\rr}\zeta(\theta;\cdot)$
is admitted to be set-valued at the magnetization rate $\rr=0$.

Further, we test \eqref{ED-5+disc} by $\mu_0(u_\etau^k{-}u_\etau^{k-1})/\tau$. By
the analog of the calculus \eqref{calculus-PM4} and exploiting convexity of the
functional $u\mapsto\int_{\R^d}\frac12\mu_0|\Nabla u|^2\,\d x$, we obtain the
inequality
\begin{align}
\nonumber&\!\!\!\!\!\!\int_{\R^d}\!\!\frac{\mu_0}{\qexp}\frac{|\Nabla u_\etau^k|^{\qexp}
-|\Nabla u_\etau^{k-1}|^2\!\!}\tau\,\d x
+\!\int_\varOmega\!\!\bigg(\mu_0
\frac{\mmk{\cdot}\hh_{\rm dem,\etau}^k\!-\mmkk{\cdot}\hh_{\rm dem,\etau}^{k-1}\!\!}\tau
\\[-.1em]\nonumber&\hspace{10em}
-\mu_0\rrk{\cdot}\hh_{\rm dem,\etau}^k
+\!\mu_0{\rm skw}(\hh_{\rm dem,\etau}^k{\otimes}\mmk){:}\Nabla\vvk\!
\\[-.3em]&\hspace{10em}
+\mu_0\vvk{\cdot}(\Nabla\hh_{\rm dem,\etau}^k)^\top\mmk\!
+\mu_0(\hh_{\rm dem,\etau}^k{\cdot}\mmk){\rm div}\,\vvk\!\bigg)\,\d x\le0\,.
\label{calculus-PM4-disc}\end{align}

The discrete analog of \eqref{formula4-} exploits the assumption of convexity
of $\FF$ and that $\varphi_\Ee'=\FF_\Ee'$ and results to the inequality
\begin{align}
\varphi_\Ee'(\Eek,\mmk){:}\frac{\Eek{-}\Eekk\!\!}\tau
+\FF_\mm'(\Eek,\mmk){\cdot}\frac{\mmk{-}\mmkk\!\!}\tau\ge
\frac{\FF(\Eek,\mmk)-\FF(\Eekk,\mmkk)\!\!}\tau\,.
\end{align}

Summing the above estimates for $k=1,...,l$ and the convexity of $\FF$, we
see the discrete analog of \eqref{energy+} as an inequality 
\begin{align}\nonumber
&\int_\varOmega\frac{\varrho}2|\vv_\etau^l|^2+\FF(\Ee_\etau^l,\mm_\etau^l)
+\frac{\kappa\mu_0}{\qexp}|\Nabla\mm_\etau^l|^{\qexp}
\,\d x+\int_{\R^d}\!\frac{\mu_0}{\qexp}|\nabla u_\etau^l|^{\qexp}\,\d x
\\[-.4em]\nonumber
&\hspace*{.1em}+\tau\!\sum_{k=1}^l\int_\varOmega\!
\GM(\theta_\etau^{k-1})|\RR_\etau^k|^2\!
+\DD|\EE(\vvk)|^2\!+\NU|\Nabla\EE(\vvk)|^p\!
+\partial_{\rr}\zeta\big(\theta_\etau^{k-1};\rrk\big){\cdot}\rrk\!
+\varkappa|\Nabla\RR_\etau^k|^2\,\d x
\\[-.1em]\nonumber
&\hspace*{0em}\le\int_\varOmega\mu_0\hh_{\text{\rm geo},\tau}^l{\cdot}\mm_\etau^l
+\frac{\varrho}2|\vv_0|^2+\FF(\Ee_0,\mm_0)
+\frac{\kappa\mu_0}{\qexp}|\Nabla\mm_0|^{\qexp}
-\mu_0\hh_{\text{\rm geo},\tau}^0{\cdot}\mm_0\,\d x
\\[-.4em]\nonumber
&\hspace*{.1em}
+\int_{\R^d}\!\frac{\mu_0}2|\nabla u_0|^2\,\d x
+\tau\sum_{k=1}^l\int_\varOmega\!\bigg(\varrho\ff_\etau^k{\cdot}\vvk
+\frac{\hh_{\text{\rm geo},\tau}^k{-}\hh_{\text{\rm geo},\tau}^{k-1}\!}\tau{\cdot}\mmkk\!
-\GG_\mm'(\mmk){\cdot}\frac{\mmk{-}\mmkk\!\!\!}\tau
\\[-.2em]&\hspace*{2.1em}
+\thetak\CCC_\eps'(\mmk){\cdot}\Big(\frac{\mmk{-}\mmkk\!\!\!}\tau
 +(\vvk{\cdot}\Nabla)\mmk\Big)
+\big(\thetak\CCC_\eps(\mmk){+}\phi(\thetak)\big){\rm div}\,\vvk\!\bigg)
\d x\,,\!
\label{mag-mech-engr-disc}\end{align}
where $u_0\in H^1(\R^d)$ denotes the weak solution to \eqref{Max+} with
$\mm=\mm_0$. For $\hh_{\text{\rm geo},\tau}^0$, we consider an (arbitrary)
extension of $\hh_\text{\rm geo}$ from $I=[0,T]$ to $[-\tau,T]$ by continuity
so that $\hh_{\text{\rm geo},\tau}^0\to \hh_\text{\rm geo}(0)$ for $\tau\to0$ in
$H^1(\varOmega;\R^d)$. Adding \eqref{ED-6+disc} tested by 1 summed for
$k=1,...,l$, we see cancellation of the adiabatic terms
$\thetak\CCC_\eps'(\mmk){\cdot}((\mmk{-}\mmkk)/\tau
 +(\vvk{\cdot}\Nabla)\mmk)$ and $(\thetak
\CCC_\eps(\mmk){+}\phi(\thetak)){\rm div}\,\vvk$
and an ``$\eps$-cancelation'' of the dissipative terms.
We thus obtain an analog of the total energy balance \eqref{energy-mag} except
that an $\eps$-part of the dissipation still remains: 
\begin{align}\nonumber
&\int_\varOmega\frac{\varrho}2|\vv_\etau^l|^2\!+\FF(\Ee_\etau^l,\mm_\etau^l)
+\frac{\kappa\mu_0}{\qexp}|\Nabla\mm_\etau^l|^{\qexp}\!
+\W_\etau^l\,\d x+\int_{\R^d}\!\frac{\mu_0}2|\nabla u_\etau^l|^2\,\d x
\\[-.1em]\nonumber&\hspace*{0em}
+\eps\tau\!\sum_{k=1}^l\int_\varOmega\!\GM(\theta_\etau^{k-1})|\RR_\etau^k|^2\!+
\DD|\EE(\vvk)|^2\!+\NU|\Nabla\EE(\vvk)|^p\!
+\partial_{\rr}\zeta\big(\theta_\etau^{k-1};
\rrk\big){\cdot}\rrk\!+\varkappa|\Nabla\RR_\etau^k|^2\,\d x
\\[-.2em]\nonumber
&\hspace*{.0em}
\le\int_\varOmega\!\mu_0\hh_{\rm ext\tau}^l{\cdot}\mm_\etau^l\!\!
+\frac{\varrho}2|\vv_0|^2\!+\FF(\Ee_0,\mm_0)
+\frac{\kappa\mu_0}{\qexp}|\Nabla\mm_0|^{\qexp}\!
-\mu_0\hh_{\text{\rm geo},\tau}^0{\cdot}\mm_0+\W_0\,\d x
\\[-.4em]
&\hspace*{0em}
+\!\int_{\R^d}\!\frac{\mu_0}2|\nabla u_0|^2\,\d x
+\tau\sum_{k=1}^l\int_\varOmega\!\!\varrho\ff_\etau^k{\cdot}\vvk\!
+\frac{\hh_{\text{\rm geo},\tau}^k{-}\hh_{\text{\rm geo},\tau}^{k-1}\!\!}\tau{\cdot}\mmkk\!
-\GG_\eps'(\mmk){\cdot}\frac{\!\mmk{-}\mmkk\!\!}\tau\,\d x.\!
\label{tot-engr-disc}\end{align}
The regularization \eqref{regul} makes $\GG_\eps'(\mmk)$ bounded for fixed
$\eps>0$ while the time-difference of $\mm$ can be estimated if written as
\begin{align}\label{pdt-m}
\frac{\!\mmk{-}\mmkk\!\!}\tau=\rrk-(\vvk{\cdot}\Nabla)\mmk+{\rm skw}(\Nabla\vvk)\mmk\,.
\end{align}

Then, by the discrete Gronwall inequality and by the coercivity of $\varphi$, 
we obtain the a-priori estimates
\begin{subequations}\label{est-disc}\begin{align}\label{est-disc1}
&\|\vv_\etau^{}\|_{L^\infty(I;L^2(\varOmega;\R^d))\,\cap L^p(I;W^{2,p}(\varOmega;\R^d))}^{}
\le C_\eps,
  \\\label{est-disc5}
&\|\Ee_\etau^{}\|_{L^\infty(I;L^2(\varOmega;\R_{\rm sym}^{d\times d}))}^{}\le C_\eps\,,
   \\\label{est-disc4}
   &\|\mm_\etau^{}\|_{L^\infty(I; H^1(\varOmega;\R^d))}^{}\le C_\eps\,,
   \\\label{est-disc7}
   &\|u_\etau^{}\|_{L^\infty(I;H^1(\R^d)) }^{}\le C_\eps\,,
   \\\label{est-disc6}
   &\|\RR_\etau^{}\|_{L^2(I;H^1(\varOmega;\R_{\rm dev}^{d\times d})) }^{}\le C_\eps\,,
   \\[-.4em]\label{est-disc8}
   &\|\rr_\etau^{}\|_{L^{\rexp}(I{\times}\varOmega;\R^d) }^{}\le C_\eps
     \ \ \ \text{ and }\ \
     \Big\|\pdt{\mm_\etau^{}}+(\overlinevvtau{\cdot}\Nabla)\overlinemmtau\Big\|_{L^{\rexp}(I{\times}\varOmega;\R^d)}^{}\le C_\eps\,,
  \\[-.2em]\label{est-disc-w}
   &\|\W_\etau^{}\|_{L^\infty(I;L^1(\varOmega))}^{}\le C_\eps\ \ \text{ and }\ \ 
   \|\theta_\etau^{}\|_{L^\infty(I;L^1(\varOmega))}^{}\le C_\eps\,.
\intertext{The later estimate in \eqref{est-disc8} follows from
$\pdt{}\mm_\etau^{}+(\overlinevvtau{\cdot}\Nabla)\overlinemmtau
=\overline\rr_\etau^{}+{\rm skw}(\Nabla\overlinevvtau)\overlinemmtau$
since ${\rm skw}(\Nabla\overlinevvtau)\overlinemmtau\in
L^p(I;L^{\qexp^*}(\varOmega;\R^d))\subset L^{\rexp}(I{\times}\varOmega;\R^d)$
because $\rexp\le\min(p,\qexp^*)$ is assumed. Moreover, from
\eqref{pdt-m} and the estimates (\ref{est-disc}a,c,f), we also have}
\label{est-disc9}
   &
  \Big\|\pdt{\mm_\etau^{}}\Big\|_{L^\rexp(I{\times}\varOmega;\R^d)}^{}\le C_\eps\,.
 \end{align}\end{subequations}

Having now estimated the dissipation heat source on the right-hand side
of the discrete heat equation \eqref{ED-6+d}, 
as the next step we can use the $L^1$-technique to estimate of temperature
gradient developed by Boccardo and Gallou{\"e}t \cite{BocGal89NEPE}
exploiting sophistically Gagliardo-Nirenberg inequality, cf.\
\cite[Prop.\,8.2.1]{KruRou19MMCM}. The essence is to test the heat-transfer
equation \eqref{ED-6+d} by a smoothened Heaviside function,
say ${\mathfrak h}(\theta)=1-(1{+}\theta)^{-\epsilon}$
for $\epsilon>0$, as suggested in \cite{FeiMal06NSET}. The modification in
comparison with the usual ``heat operator'' in the form
$\frac{\partial}{\partial t}\gamma(\theta)-{\rm div}(\KK(\theta)\Nabla\theta)$
and with an $L^1$-right-hand side consists in that $\DT\W$ contains one more
term, namely ${\rm div}(\vv\gamma(\theta))$. In the discrete form, this
additional convective term is
${\rm div}(\overlinevvtau\gamma(\overline\theta_\etau))$
and, when tested by ${\mathfrak h}(\overline\theta_\etau)$, we can estimate
it ``on the right-hand side'' as
\begin{align}\nonumber
 -\!\int_\varOmega\!{\rm div}(\overlinevvtau\gamma(\overline\theta_\etau))
  {\mathfrak h}(\overline\theta_\etau)\,\d x
  &=\int_\varOmega\overlinevvtau{\cdot}\gamma(\overline\theta_\etau)
  {\mathfrak h}'(\overline\theta_\etau)\Nabla\overline\theta_\etau\,\d x
  \\&\le\int_\varOmega\frac1\epsilon|\overlinevvtau|^2
  \gamma(\overline\theta_\etau)^2{\mathfrak h}'(\overline\theta_\etau)\,\d x
  +\epsilon\int_\varOmega\!{\mathfrak h}'(\overline\theta_\etau)|\Nabla\overline\theta_\etau|^2\,\d x\,.
\label{gamma-calculus}\end{align}
As $\gamma(\theta)=\mathscr{O}(\theta)$ while
${\mathfrak h}'(\theta)=\mathscr{O}(1/\theta)$, we have
$[\gamma^2{\mathfrak h}'](\theta)=\mathscr{O}(\theta)$ so that
$\gamma(\overline\theta_\etau)^2{\mathfrak h}'(\overline\theta_\etau)$ is
bounded in $L^\infty(I;L^1(\varOmega))$ while $|\overlinevvtau|^2$
is surely bounded in $L_{\rm w*}^{p/2}(I;L^\infty(\varOmega))$, so that 
the integrand in the penultimate integral in \eqref{gamma-calculus} is bounded
in $L^{p/2}(I;L^1(\varOmega))$.
For the last integral, this is exactly fitted with the estimation in the
$L^1$-theory, and for $\epsilon>0$ sufficiently small can be absorbed in
the respective estimation, cf.\ \cite[Formula (8.2.17)]{KruRou19MMCM}.
Using also the already obtained estimate \eqref{est-disc-w}, we have a
``prefabricated'' estimate 
$\delta_s\|\theta_\etau\|_{L^s(I{\times}\varOmega;\R^d)}^s+\delta_\sigma\|\Nabla\theta_\etau\|_{L^\sigma(I{\times}\varOmega;\R^d)}^\sigma\le
1+\|{\rm rhs}_\etau\|_{L^1(I{\times}\varOmega)}$ for some $\delta_s>0$ 
at disposal with ``rhs$_\etau$'' abbreviating the right-hand side of
\eqref{ED-6+d} and with $s<1+2/d$ and $\sigma<(d{+}2)/(d{+}1)$, cf.\
\cite[Prop.\,8.2.1]{KruRou19MMCM}. We add this estimate to the discrete
mechanical energy balance \eqref{mag-mech-engr-disc} for $l=T/\tau$
multiplied by the factor 2. Assuming (witout loss of generality) $\FF\ge0$
and then forgetting the first (nonnegative) integral in
\eqref{mag-mech-engr-disc}, we thus obtain the estimate
\begin{align}\nonumber
&\hspace*{.0em}\int_0^T\!\!\!\int_\varOmega\!\Big(\!
\GM(\underline\theta_\etau)|\overlineRRtau|^2\!
+\DD|\EE(\overlinevvtau)|^2\!+\NU|\Nabla\EE(\overlinevvtau)|^p\!
\\[-.6em]\nonumber&\hspace*{8em}
+\partial_{\rr}\zeta\big(\underline\theta_\etau;\overline\rr_\etau\big){\cdot}
\overline\rr_\etau\!
+\varkappa|\Nabla\overlineRRtau|^2\!
+\delta_s|\theta_\etau|^s+\delta_\sigma|\Nabla\theta_\etau|^\sigma
\Big)\,\d x\d t
\\[-.2em]\nonumber&\hspace*{.2em}
\le\!\int_\varOmega\!2\mu_0\hh_{\text{\rm geo},\tau}(T){\cdot}\mm_\etau(T)
+{\varrho}|\vv_0|^2\!+2\FF(\Ee_0,\mm_0)
+\kappa\mu_0|\Nabla\mm_0|^{\qexp}-2\mu_0\hh_{\text{\rm geo},\tau}^0{\cdot}\mm_0\,\d x
\\[-.2em]\nonumber
&\hspace*{.2em}
+\!\int_{\R^d}\!{\mu_0}|\nabla u_0|^2\,\d x
+2\int_0^T\!\!\!\int_\varOmega\!\varrho\overline\ff_\etau{\cdot}\overlinevvtau\!
+\frac{\partial\hh_{\text{\rm geo},\tau}}{\partial t}{\cdot}\underline\mm_\etau
-\GG_\eps'(\overlinemmtau){\cdot}\frac{\partial\mm_\etau\!}{\partial t}
\\[-.2em]
&\hspace*{1.5em}
+\big|\overline\theta_\etau\CCC_\eps'(\overlinemmtau)\big|
\Big|\frac{\partial\mm_\etau\!}{\partial t}+(\overlinevvtau{\cdot}\Nabla)\overlinemmtau\Big|
+\big|\big(\overline\theta_\etau
\CCC_\eps(\overlinemmtau){+}\phi(\overline\theta_\etau)\big){\rm div}\,\overlinevvtau\big|
\bigg)\,\d x\d t\,.
\label{mag-mech-engr-disc+}\end{align}
The term $\GG_\eps'(\overlinemmtau){\cdot}\frac{\partial}{\partial t}\mm_\etau$
can be estimated by using \eqref{est-disc9} and the growth assumption
\eqref{ass:1-}. The penultimate term in \eqref{mag-mech-engr-disc+} can be
estimated by meas of the growth \eqref{ass:1--} of $\CCC'$ 
and the latter estimate in \eqref{est-disc8}, 
we have by the H\"older and the Young inequalities the estimate 
\begin{align}\nonumber
&\hspace*{.0em}\int_0^T\!\!\!\int_\varOmega\!
\big|\overline\theta_\etau\CCC_\eps'(\overlinemmtau)\big|
\Big|\frac{\partial\mm_\etau\!}{\partial t}+(\overlinevvtau{\cdot}\Nabla)\overlinemmtau\Big|\,\d x\d t
\\[-.4em]&\qquad\le
\frac{\delta_s}3\big\|\overline\theta_\etau\big\|_{L^s(I{\times}\varOmega)}^s+C
\big(1{+}\|\overlinemmtau\|_{L^{\qexp^*}(I{\times}\varOmega;\R^d)}^{s'}\big)
\Big\|\frac{\partial\mm_\etau\!}{\partial t}{+}(\overlinevvtau{\cdot}\Nabla)\overlinemmtau\Big\|_{L^{\rexp}(I{\times}\varOmega;\R^d)}^{s'}
\label{adiabat-est-1}\end{align}
provided $1/s+1/\qexp^*+1/\rexp\le1$, which needs a sufficiently big  $r$ as
specified in \eqref{ass:4}. Actually, the growth of $\GG_\eps'$ and
$\CCC_\eps'$ with $\GG_\eps$ and $\CCC_\eps$ from \eqref{regul} is even better
but, in Step~5, we will
use that this estimate is even uniform with respect to $\varepsilon>0$.
The last term in \eqref{mag-mech-engr-disc+} can be estimated by using
the growth assumption (\ref{ass}c,d) on $\CCC$ and on $\phi$ as
\begin{align}\nonumber
&\hspace*{.0em}\int_0^T\!\!\!\int_\varOmega\!
\big|\big(\overline\theta_\etau
\CCC_\eps(\overlinemmtau){+}\phi(\overline\theta_\etau)\big)
{\rm div}\,\overlinevvtau\big|\,\d x\d t
\\[-.4em]&\qquad\quad\le
\frac{\delta_s}3\|\overline\theta_\etau\big\|_{L^s(I{\times}\varOmega)}^s+C
\big(1{+}\|\overlinemmtau\|_{L^\infty(I;L^{\qexp^*}_{^{^{}}}(\varOmega;\R^d))}^{s'}\big)
\|\Nabla\overlinevvtau\|_{L^p(I;L^\infty(\varOmega;\R^{d\times d}))}^{s'}\,,
\label{adiabat-est-2}\end{align}
which holds for $1/s+1/p\le1$ (which is satisfied if $p>d$)
and for $1/s+1/\qexp^*\le1$ (which is satisfied if $s\ge{\qexp^*}'$).
Merging \eqref{adiabat-est-1}--\eqref{adiabat-est-2} with 
\eqref{mag-mech-engr-disc+}, we obtain the estimates:
\begin{subequations}\label{est-disc+}\begin{align}\label{est-disc+theta}
&\|\overline\theta_\etau\|_{L^\sigma(I;W^{1,\sigma}(\varOmega))
 \,\cap\,L^s(I{\times}\varOmega)}^{}\le C_{\eps,s,\sigma}\,,
\\[-.4em]\label{est-disc+vartheta}
&\|\overline\W_\etau\|_{L^\sigma(I;W^{1,\sigma}(\varOmega))
  \,\cap\,L^s(I{\times}\varOmega)}^{}\le C_{\eps,s,\sigma}
\ \ \text{ with $\ \ 1\le s<1+\frac2d$\,,\
$\ 1\le\sigma<\frac{d+2}{d+1}$.}
\end{align}\end{subequations}
The $L^\sigma$-estimate \eqref{est-disc+vartheta}of $\nabla\overline\W_\etau$ 
can be read from \eqref{est-disc+theta}
due to $\nabla\overline\W_\etau=\nabla\gamma(\overline\theta_\etau)
=\gamma'(\overline\theta_\etau)\nabla\overline\theta_\etau$.
The $L^s$-estimate of $\overline\W_\etau$ 
is to be read by the Gagliardo-Nirenberg interpolation
of the first estimate in \eqref{est-disc+vartheta} with  \eqref{est-disc-w}.

In addition, by comparison from \eqref{ED-4+d}, we obtain also
\begin{align}
&\label{est-Delta}
\big\|\Delta\overlinemmtau\big\|_{L^{\rexp'}(I{\times}\varOmega;\R^d)}^{}\le C_\eps\,.
 \end{align}
More in detail, using also the assumption \eqref{ass:2} and the bound
\eqref{est-disc8}, we have $\{\partial_{\ZJ\mm}\zeta(\underline\theta_\etau,\overline\rr_\etau)\}_{\tau>0}$ bounded in $L^{\rexp'}(I{\times}\varOmega;\R^d)$.
Also $\{\Nabla u|_{I{\times}\varOmega}^{}\}_{\tau>0}$ and
$\{\FF_\mm'(\overlineEetau,\overlinemmtau)+[\varpi_\eps]_\mm'(\overlinemmtau,\overline\theta_\etau)\}_{\tau>0}$ are bounded in
$L^{\rexp'}(I{\times}\varOmega;\R^d)$ due to \eqref{ass:1-};
here we used the growth assumption \eqref{ass:1--} on
$\CC'$ together with the estimates \eqref{est-disc4}
and \eqref{est-disc+theta}, so that
\begin{align}\nonumber
\int_0^T\!\!\!\int_\varOmega\big|\overline\theta_\etau\CC_\eps(\overlinemmtau)
\big|^{\rexp'}\,\d x\d t
&\le C\!\int_0^T\!\!\!\int_\varOmega1{+}\overline\theta_\etau^{\,\rexp'}\big|\overlinemmtau\big|^{\rexp'}\,\d x\d t
\\&\le C'\big(1+\|\overline\theta_\etau\|_{L^{1+2/d-\epsilon}(I{\times}\varOmega)}^{\rexp'}
\|\overlinemmtau\|_{L^{\qexp^*}(I{\times}\varOmega)}^{\rexp'}\big)
\label{est-of-theta.omega}\end{align}
provided $d/(d{+}2)+1/\qexp^*<1/\rexp'$, which gives the restriction
$\rexp>\frac{\qexp^*d{+}\qexp^*2}{\qexp^*d{+}\qexp^*2{-}\qexp^*d{-}d{-}2}$,
i.e.\ $\rexp>30/7$ for $d=3$ or $\rexp>2$ for $d=2$ as assumed in
\eqref{ass:4}.

\def\expSIGMA{\sigma}
An important attribute of the model is that the convective transport of
variables via the velocity field
$\vv\in L_{\rm w*}^1(I;W^{1,\infty}(\varOmega;\R^d))$ or, in the discrete variant
by $\overline\vv_\tau$ bounded in $L^1(I;W^{1,\infty}(\varOmega;\R^d))$,
qualitatively well copies regularity properties of the initial conditions.
We use this phenomenon particularly for the Zaremba-Jaumann time difference
of the elastic strain $\Ee=[E_{ij}]$.
Let us consider a general tensor-valued source ${\bm F}\in
L^2(I{\times}\varOmega;\R^{d\times d})$ for $\ZJ\Ee={\bm F}$, i.e.\ in the difference variant
\begin{align}
\frac{{\bm E}_\tau^k-{\bm E}_\tau^{k-1}\!\!}\tau
  +B_\text{\sc zj}(\vv_\tau^k,{\bm E}_\tau^k)={\bm F}_\tau^k\,.
  \label{ZJ-evol-abstract}\end{align}
For $\expSIGMA>1$, we use the following calculus exploiting the Green formula with
the boundary condition $\vv{\cdot}\nn=0$: 
\begin{align}\nonumber
  \int_{\varOmega}(\vv{\cdot}\nabla z)|z|^{\expSIGMA-2}z\,\d x
  &=\int_{\varOmega}\!(1{-}\expSIGMA)|z|^{\expSIGMA-2}z(\vv{\cdot}\nabla z)-({\rm div}\,\vv)|z|^\expSIGMA\d x
\\[-.4em]&\hspace{9em}+\int_{\varGamma}{}\!|z|^\expSIGMA(\vv{\cdot}\nn)\,\d S
=-\frac1\expSIGMA\int_{\varOmega}({\rm div}\,\vv)|z|^\expSIGMA\d x\,.
\label{transport}\end{align}
We test \eqref{ZJ-evol-abstract} by $|{\bm E}_\tau^k|^{\expSIGMA-2}{\bm E}_\tau^k$,
which gives
\begin{align}
  \frac1\expSIGMA\int_{\varOmega}\!\frac{|{\bm E}_\tau^k|^\expSIGMA\!-|{\bm E}_\tau^{k-1}|^\expSIGMA\!}\tau
  \,\d x\le\int_{\varOmega}\!\frac{{\rm div}\,\vv_\tau^k}\expSIGMA|{\bm E}_\tau^k|^\expSIGMA\!
  +2|{\rm skw}(\nabla\vv_\tau^k)|\,|{\bm E}_\tau^k|^\expSIGMA\!
  +|{\bm E}_\tau^k|^{\expSIGMA-2}{\bm E}_\tau^k{:}{\bm F}_\tau^k\,\d x,
\label{transport+}\end{align}
where we used \eqref{transport} for each component $z=E_{ij}$. From
\eqref{transport+}, by the Young and the discrete Gronwall inequalities,
we obtain the estimate
\begin{align}
  \|{\bm E}_\tau^k\|_{L^\expSIGMA(\varOmega;\R^{d\times d})}^\expSIGMA\le
 C{\rm e}^{1+2k\tau\max_{l=1,...,k}\|\nabla\vv_\tau^k\|_{L^\infty(\varOmega;\R^{d\times d})}}\bigg(\|{\bm E}_\tau^0\|_{L^\expSIGMA(\varOmega;\R^{d\times d})}^\expSIGMA\!+
  \tau\sum_{l=1}^k\|{\bm F}_\tau^l\|_{L^{\expSIGMA'}(\varOmega;\R^{d\times d})}^{\expSIGMA'}\bigg)\,
\label{Gronwall}\end{align}
for some $C$ and for $\tau\le1/
((2{+}4\expSIGMA)\max_{l=1,...,k}\|\nabla\vv_\tau^k\|_{L^\infty(\varOmega;\R^{d\times d})}+2\expSIGMA)$.

Moreover, we can also test \eqref{ZJ-evol-abstract} by the $q$-Laplacian
$-{\rm div}(|\nabla{\bm E}_\tau^k|^{q-2}\nabla{\bm E}_\tau^k)$. More specifically,
we can apply the $\nabla$-operator to \eqref{ZJ-evol-abstract}
and test it by $|\nabla{\bm E}_\tau^k|^{q-2}\nabla{\bm E}_\tau^k$.
Instead of \eqref{transport}, we use the calculus 
\begin{align}\nonumber
&\int_{\varOmega}\nabla(\vv{\cdot}\nabla z){\cdot}|\nabla z|^{q-2}\nabla z\,\d x
=\int_{\varOmega}|\nabla z|^{q-2}\nabla\vv{:}(\nabla z{\otimes}\nabla z)
+(\vv{\cdot}\nabla^2z){\cdot}|\nabla z|^{q-2}\nabla z\,\d x
\\[-.4em]&\nonumber\quad=
-\int_{\varOmega}\!|\nabla z|^{q-2}\nabla\vv{:}(\nabla z{\otimes}\nabla z)
+(q{-}1)|\nabla z|^{q-2}\nabla z{\cdot}(\vv{\cdot}\nabla^2z)
+({\rm div}\,\vv)|\nabla z|^q\,\d x
\\[-.4em]&\hspace{4.7em}
+\int_{\varGamma}{}\!|\nabla z|^q(\vv{\cdot}\nn)\,\d S
=\int_{\varOmega}|\nabla z|^{q-2}\nabla\vv{:}(\nabla z{\otimes}\nabla z)
-\frac1q({\rm div}\,\vv)|\nabla z|^q\,\d x\,.
\label{transport-}\end{align}
In the tensorial situation \eqref{ZJ-evol-abstract}, we use it again for
$z=E_{ij}$ and then we use also 
\begin{align}\nonumber
&\nabla\Big({\rm skw}(\nabla\vv_\tau^k){\bm E}_\tau^k
-{\bm E}_\tau^k{\rm skw}(\nabla\vv_\tau^k)\Big)
\Vdots|\nabla{\bm E}_\tau^k|^{q-2}\nabla{\bm E}_\tau^k
\\[-.1em]&\hspace{14em}
\le2|\nabla\vv_\tau^k|\,|\nabla{\bm E}_\tau^k|^q
+2|\nabla^2\vv_\tau^k|\,|{\bm E}_\tau^k|\,|\nabla{\bm E}_\tau^k|^{q-1}.
\label{transport--}\end{align}
Instead of \eqref{transport+}, this gives
\begin{align}\nonumber
\frac1q\int_{\varOmega}\!
\frac{|\nabla {\bm E}_\tau^k|^q\!-|\nabla {\bm E}_\tau^{k-1}|^q\!}\tau
\,\d x&\le\int_{\varOmega}\!\bigg(\frac{{\rm div}\,\vv_\tau^k}q|\nabla{\bm E}_\tau^k|^q
+\Big(2+\frac1q\Big)|\nabla\vv_\tau^k|
\,|\nabla{\bm E}_\tau^k|^q
\\[-.7em]\nonumber&\qquad
+|\nabla{\bm E}_\tau^k|^{q-2}\nabla{\bm E}_\tau^k{\Vdots}\nabla{\bm F}_\tau^k
+2|\nabla^2\vv_\tau^k|\,|{\bm E}_\tau^k|\,|\nabla{\bm E}_\tau^k|^{q-1}
\bigg)\,\d x
\\[-.3em]&\hspace{-9em}\nonumber\le
C+C\big(1{+}\|\nabla\vv_\tau^k\|_{L^\infty(\varOmega;\R^{d\times d})}^{}\big)
\|\nabla{\bm E}_\tau^k\|_{L^q(\varOmega;\R^{d\times d\times d})}^q
\\[-.1em]
&\hspace{-4em}
+\|\nabla{\bm F}_\tau^k\|_{L^q(\varOmega;\R^{d\times d\times d})}^q+
\|\nabla^2\vv_\tau^k\|_{L^p(\varOmega;\R^{d\times d\times d})}^p
+\|{\bm E}_\tau^k\|_{L^\expSIGMA(\varOmega;\R^{d\times d\times d})}^\expSIGMA
\,.
\label{transport+++}\end{align}
for some $C$ sufficiently large, provided $1/p+1/\expSIGMA+1/q'\le1$ with $p$
from \eqref{est-disc1}. By a discrete Gronwall inequality like
\eqref{Gronwall}, provided $\tau$ is sufficiently small, we obtain
\begin{align}
  \|\nabla{\bm E}_\tau^k\|_{L^q(\varOmega;\R^{d\times d\times d})}^q\le  C
\bigg(\|\nabla{\bm E}_\tau^0\|_{L^q(\varOmega;\R^{d\times d\times d})}^q\!+
\tau\sum_{l=1}^k\|\nabla{\bm F}_\tau^l\|_{L^{q}(\varOmega;\R^{d\times d\times d})}^{q}\bigg)\,
\label{Gronwall+}\end{align}
with some $C$ depending on
$\max_{l=1,...,k}\|\nabla\vv_\tau^k\|_{L^\infty(\varOmega;\R^{d\times d})}$ and on 
$\tau\sum_{l=1}^k\|\nabla^2\vv_\tau^k\|_{L^p(\varOmega)}^p\!+\|{\bm E}_\tau^k\|_{L^\expSIGMA(\varOmega)}^\expSIGMA$.

Considering $\vv_\tau^k=\vvk$ and ${\bm F}_\tau^k=\EE(\vvk)-\RR_\etau^k$
and using the already obtained estimates (\ref{est-disc}a,e) and
the initial condition $\Ee_0\in H^1(\varOmega;\R_{\rm dev}^{d\times d})$, 
we use the calculus \eqref{transport+}--\eqref{Gronwall}
with $\expSIGMA=2^*$ and $q=2$ for \eqref{ED-2+disc}. 
Thus we obtain
\begin{align}\label{Gronwall-3-2}
\|\nabla\overlineEetau\|_{L^\infty(I;L^2(\varOmega;\R^{d\times d\times d}))}^{}\le C\,.
\end{align}

\def\whvv{\widehat{\vv}_\tau}

\medskip\noindent{\it Step 3:
Convergence in the mechanical part for $\tau\to0$ with $\eps>0$ fixed.}
By the Banach selection principle, we obtain a subsequence
converging weakly* with respect to the topologies indicated in
\eqref{est-disc} and \eqref{est-disc+} to some limit
$(\vv_\eps,\Ee_\eps,\RR_\eps,\mm_\eps,u_\eps,\theta_\eps)$.

Moreover, we use the Aubin-Lions compact-embedding theorem generalized for
functions with measure time derivatives as in \cite[Cor.\,7.9]{Roub13NPDE}
to show the strong convergence 
\begin{subequations}\label{conv-strong}\begin{align}\label{conv-strong-nabla-v}
    &&&\Nabla\overlinevvtau\to\Nabla\vv_\eps
&&\text{ in }\ L^2(I;L^\infty(\varOmega;\R^{d\times d}))\,,
&&&&
\\\label{conv-strong-E}
&&&\overlineEetau\to\Ee_\eps
&&\text{ in }\ L^{1/\epsilon}(I;L^{2^*-\epsilon}(\varOmega;\R_{\rm sym}^{d\times d}))\,,
&&&&
\\&&&\label{conv-strong-a}
    \overlinemmtau\to\mm_\eps&&\text{ in }\
    L^{1/\epsilon}(I;L^{\qexp^*-\epsilon}(\varOmega;\R^d))\,,
\\&&&\label{conv-strong-+++}
\varphi(\overlineEetau,\overlinemmtau)\to
\varphi(\Ee_\eps,\mm_\eps)&&\text{ in }\ L^{1/\epsilon}(I;L^1(\varOmega))\,,
&&&&
\\&&&\label{conv-strong-+}
\varphi_\mm'(\overlineEetau,\overlinemmtau)\to\varphi_\mm'(\Ee_\eps,\mm_\eps)\!\!\!\!&&\text{ in }\ L^{1/\epsilon}(I;L^{\rexp'}(\varOmega;\R^d))\,,
&&&&
\\&&&\label{conv-strong-omega-1}\CCC_\eps(\overlinemmtau)\to
\CCC_\eps(\mm_\eps)&&\text{ in }\ L^{1/\epsilon}(I;L^{\qexp^*/2-\epsilon}(\varOmega))\,,
\\&&&\label{conv-strong-omega-1+}\CCC_\eps'(\overlinemmtau)\to
\CCC_\eps'(\mm_\eps)&&\text{ in }\
L^{1/\epsilon}(I;L^{\qexp^*-\epsilon}(\varOmega;\R^d))
\end{align}\end{subequations}
for any $0<\epsilon\le\qexp^*{-}1$;
here for (\ref{conv-strong}d--f) we used the growth assumptions (\ref{ass}a,c).
For (\ref{conv-strong}a--c), we need some information about time derivatives
to be able to use the mentioned Aubin-Lions theorem generalized
for piece-wise
constant functions in time. For \eqref{conv-strong-nabla-v}, a boundedness
of the sequence $\{\pdt{}\vv_\etau\}_{\tau>0}$ in $L^1(I;W^{2,p}(\varOmega;\R^d)^*)$
can be seen by comparison from \eqref{ED-1+d}. For the convergence
\eqref{conv-strong-E}, we can see some information by comparison from
\eqref{ED-2+d}, from which we can see that 
$\{\frac{\partial}{\partial t}\Ee_\etau\}_{\tau>0}$ is bounded in
$L^2(I;H^1(\varOmega;\R_{\rm sym}^{d\times d})^*)$. Eventually, for
\eqref{conv-strong-a}, an estimate of $\pdt{}\mm_\etau$
is directly at disposal in the a-priori estimate
\eqref{est-disc9}. The strong convergence in
(\ref{conv-strong}d--f) is then simply by continuity of the Nemytski\u{\i}
operator induced by $\varphi$, $\varphi_\Ee'$, and $\varphi_\mm'$
even without having any information about the time derivative.

As $\Nabla\mm$ occurs nonlinearly in the capillarity stress and also in the
weak formulation \eqref{VI} multiplied by $\Delta\mm$, we need to prove also a
strong convergence $\Nabla\overlinemmtau\to\Nabla\mm_\eps$.
To prove it, we take a sequence $\{\widetilde\mm_\tau\}_{\tau>0}^{}$
piecewise constant in time with respect to the partition with the time
step $\tau$  and converging strongly towards $\mm_\eps$ for $\tau\to0$.
Then, using the variational inequality arising from the inclusion
\eqref{ED-4+d} tested by $\overlinemmtau{-}\widetilde\mm_\tau$,
we can see that, written in terms of the interpolants, that
it holds
\begin{align} \nonumber
&\int_0^T\!\!\!\int_\varOmega 
\kappa|\Nabla(\overlinemmtau{-}\widetilde\mm_\tau)|^{\qexp}\,\d x\d t
=\int_0^T\!\!\!\int_\varOmega
\bigg(\frac{\varphi_\mm'(\overlineEetau,\overlinemmtau){+}
\overline\theta_\etau\CC_\eps'(\overlinemmtau)
\!}{\mu_0}
-\partial_{\rr}^{}\zeta\big(\underline\theta_\etau;
\overline\rr_\etau\big)
\\[-.3em]&\hspace{5em}
+\overline\hh_{\text{\rm geo},\tau}-\Nabla\overline u_\etau
\bigg){\cdot}(\overlinemmtau{-}\widetilde\mm_\tau)-\kappa
\Nabla\widetilde\mm_\etau{:}
\Nabla(\overlinemmtau{-}\widetilde\mm_\tau)\,\d x\d t\to0\,.
\label{strong-nabla-mm}
\end{align}
In fact, here we mean a suitable selection from the set
$\partial_{\rr}^{}\zeta(\underline\theta_\etau;\overline\rr_\etau)$.
Here we again used \eqref{ass:2} with \eqref{est-disc8} so that the set
$\partial_{\rr}^{}\zeta(\underline\theta_\etau;\overline\rr_\etau)$ is bounded
(uniformly with respect to $\tau$) in $L^{\rexp'}(I{\times}\varOmega;\R^d)$.
Moreover, we used growth assumptions (\ref{ass}a,c,f) and the a-priori estimates
(\ref{est-disc}c,d), \eqref{est-disc+}, and \eqref{Gronwall-3-2}, so that
$\varphi_\mm'(\overlineEetau,\overlinemmtau)$ and
$\overline\theta_\etau\CC_\eps'(\overlinemmtau)$ are bounded in
$L^{\rexp'}(I{\times}\varOmega;\R^d)$ 
while $\overlinemmtau-\widetilde\mm_\tau\to0$ strongly surely in
$L^{\qexp^*-\epsilon}(I{\times}\varOmega;\R^d)\subset
L^\rexp(I{\times}\varOmega;\R^d)$; cf.\ \eqref{conv-strong-a}.
In \eqref{strong-nabla-mm}, we used also that both
$\Nabla\overline u_\etau$ is bounded in $L^2(I{\times}\varOmega;\R^d)
\subset L^{\rexp'}(I{\times}\varOmega;\R^d)$; here we used
$1/\qexp^*+1/\rexp'<1$ which is granted if $r<\qexp^*$, cf.\ the assumption
\eqref{ass:4}. We thus obtained the desired strong convergence in 
$L^\qexp(I{\times}\varOmega;\R^{d\times d})$. By interpolation with
the estimate \eqref{est-disc4}, we obtain even
\begin{align}
&&&\Nabla\overlinemmtau\to\Nabla\mm_\eps&&\hspace*{-2em}\text{strongly in }\ 
L^{1/\epsilon}(I;L^\qexp(\varOmega;\R^{d\times d}))
\,.&&&&
\label{nabla-dam-strong}\end{align}

As $\GM(\cdot)$, $\zeta(\cdot;\ZJ\mm)$, and $\phi(\cdot)$ depend nonlinearly
on the temperature $\theta$, we need a strong convergence of
$\overline\theta_\etau$ and $\underline\theta_\etau$:
\begin{align}
\overline\theta_\etau\to\theta_\eps\ \text{ and }\ \ \underline\theta_\etau\to\theta_\eps
\ \text{ strongly in }\ L^s(I{\times}\varOmega)
\label{beta-discrete}\end{align}
with $s$ from 
\eqref{est-disc+vartheta}. This follows from the estimates on the gradient
\eqref{est-disc+vartheta}, when using the (generalized) Aubin-Lions compact
embedding theorem interpolated also with the a-priori estimates
\eqref{est-disc-w}. Then we obtain also
\begin{subequations}\label{conv-strong+}\begin{align}
&&&\label{conv-strong-+-}
\overline\theta_\etau\CC_\eps(\overlinemmtau)\to\theta_\eps\CC_\eps(\mm_\eps)
\!\!\!\!\!&&
\text{ in }\ L^s(I;L^{\qexp^*s/(\qexp^*+2s)}(\varOmega;\R^d))\,,
\\&&&\label{conv-strong-++}
\overline\theta_\etau\CC_\eps'(\overlinemmtau)\to\theta_\eps\CC_\eps'(\mm_\eps)
\!\!\!\!\!&&
\text{ in }\ L^{\qexp^*s/(\qexp^*+s)}(I{\times}\varOmega;\R^d)\,,
&&&&
\\&&&\label{conv-strong-++-}
\phi(\overline\theta_\etau)\to
\phi(\theta_\eps)&&\text{ in }\ L^s(I{\times}\varOmega;\R^d)\,
&&&&
\end{align}\end{subequations}
with $s$ from \eqref{est-disc+}. Note that, for $d=3$, $s<5/3$ and $\qexp^*=6$
so that the exponent in \eqref{conv-strong-+-} is below 15/14 but still
bigger than 1. Therefore, in addition to \eqref{conv-strong}, we now have
also 
\begin{align}
&&&\overlineSetau\to\bm{S}_{\text{\sc e},\eps}=\varphi_\Ee'(\Ee_\eps,\mm_\eps)+
\psi(\Ee_\eps,\mm_\eps,\theta_\eps)\bbI
&&\text{ in }\ L^1(I{\times}\varOmega;\R_{\rm sym}^{d\times d})\,.
&&&&
\label{conv-strong-S}
\end{align}

Using also \eqref{conv-strong-+++}, \eqref{nabla-dam-strong},
\eqref{beta-discrete}, and \eqref{conv-strong+}, we can see that the
capillarity/couple stress converges as 
\begin{align}\nonumber
\overlineSstrtau\to{\bm S}_{\text{\sc c},\eps}^{}&=
\kappa\mu_0\Big(\Nabla\mm_\eps{\otimes}\Nabla\mm_\eps-
\frac1{\qexp}|\Nabla\mm_\eps|^2\bbI\Big)-\mu_0{\rm skw}\big(\hh_{\rm drv,\eps}
{\otimes}\mm_\eps\big)
\\&\hspace*{9em}
  \text{ weakly in }\ L^{\min(s,\rexp')}(I;L^1(\varOmega;\R^{d\times d}))\,.
\label{Sstr-conv}\end{align}
The term which makes this convergence only weak is
${\rm skw}(\Delta\overlinemmtau{\otimes}\overlinemmtau)$
contained in $\overlineSstrtau$ from \eqref{ED-1+d} with $\overline\hh_{\rm drv,\etau}$ from \eqref{ED-4+d}; note that
$\Delta\overlinemmtau$ is bounded (and converges weakly)
in $L^{\rexp'}(I{\times}\varOmega;\R^d)$
while $\overlinemmtau$ converges strongly in
$L^{1/\epsilon}(I;L^{\qexp^*-\epsilon}(\varOmega;\R^d))$ and that $1/\rexp'+1/\qexp^*\le1$ because $\qexp^*\ge r$ is assumed. Furthermore, the term
$\mu_0\kappa(\Nabla\overlinemmtau{\otimes}\Nabla\overlinemmtau\!
-|\Nabla\overlinemmtau|^\qexp\bbI/\qexp)$
converges strongly in $L^{1/\epsilon}(I;L^1(\varOmega;\R_{\rm sym}^{d\times d}))$, 
cf.\ \eqref{nabla-dam-strong}, and is bounded
in $L^\infty(I;L^1(\varOmega;\R^{d\times d}))$.
The weak convergence is again due to the term
$\Delta\overlinemmtau{\cdot}\overlinemmtau$ contained in
$\overline\hh_{\rm drv,\etau}{\cdot}\overlinemmtau$.

By the uniform monotonicity of the operators $\frac{\partial}{\partial t}$
and ${\rm div}({\rm div}(\NU|\Nabla\EE(\cdot)|^{p-2}\Nabla\EE(\cdot))-
\DD\EE(\cdot))$, we obtain the strong convergence
\begin{subequations}\label{velocity-strong}\begin{align}
&&&\overlinevvtau\to\vv_\eps&&\text{in }\ L^p(I;W^{2,p}(\varOmega;\R^d))\ \ \text{ and }&&&&
    \\&&&\vv_\etau(T)\to\vv_\eps(T)&&\text{in }\ L^2(\varOmega;\R^d)\,.
\end{align}\end{subequations}
More in detail,
we use the discrete momentum equation \eqref{ED-1+d} tested by
$\overlinevvtau{-}\whvv$ with some
$\whvv$ piecewise constant on the time-partition of the
time step $\tau$ and converging strongly to $\vv_\eps$ in
$L^p(I;W^{2,p}(\varOmega;\R^d))\cap
C_{\rm w}(I;L^2(\varOmega;\R^d))$.
After integration over time and using also \eqref{Green-for-momentum}
to avoid usage of $\Nabla^2u$, we obtain 
\begin{align}\nonumber
  &\frac\varrho2\|\vv_\etau(T){-}\whvv(T)\|_{L^2(\varOmega;\R^d)}^2+
\NU\|\Nabla\EE(\overlinevvtau{-}\whvv)\|_{L^p(I{\times}\varOmega;\R^{d\times d\times d})}^p
  \\[-.3em]&\nonumber
  \le \int_\varOmega\frac\varrho2|\vv_\etau(T){-}\whvv(T)|^2\d x
 +\int_0^T\!\!\!\int_\varOmega\!\bigg(\DD\EE(\overlinevvtau{-}\whvv){:}
  \EE(\overlinevvtau{-}\whvv)
   \\[-.4em]&\hspace*{5em}\nonumber
   +\NU\big(|\Nabla\EE(\overlinevvtau)|^{p-2}\Nabla\EE(\overlinevvtau)
  -|\Nabla\EE(\whvv)|^{p-2}\Nabla\EE(\whvv)\big)\Vdots
  \Nabla\EE(\overlinevvtau{-}\whvv)\bigg)\,\d x\d t
 \\[-.5em]&\le\nonumber
  \int_0^T\!\!\!\int_\varOmega\!\bigg(
  \varrho\Big(\overline\ff_\tau-\frac12({\rm div}\,\overlinevvtau)\,\overlinevvtau+\mu_0(\Nabla\overline\hh_{\text{\rm geo},\tau})^\top\overlinemmtau\Big)
  {\cdot}(\overlinevvtau{-}\whvv)
\\[-.1em]\nonumber&\hspace{3em}
 -\big(\overlineSetau\!+\overlineSstrtau\!-\DD\EE(\whvv)\big){:}\EE(\overlinevvtau{-}\whvv)
-\NU|\Nabla\EE(\whvv)|^{p-2}\Nabla\EE(\whvv)\Vdots
 \Nabla\EE(\overlinevvtau{-}\whvv)
\\[-.3em]\nonumber&\hspace{3em}
-\mu_0(\Nabla\overline u_\etau{\cdot}\overlinemmtau){\rm div}(\overlinevvtau{-}\whvv)-((\overlinevvtau{-}\whvv){\otimes}\Nabla\overline u_\etau){:}\Nabla\overlinemmtau\!\bigg)
\,\d x\d t
\\[-.4em]&\hspace{18em}
  -\!\int_\varOmega\varrho\whvv(T){\cdot}(\vv_\etau(T){-}\whvv(T))\,\d x\to0\,.
\label{strong-hyper}\end{align}
Using (\ref{conv-strong}a,c,f)
and \eqref{Sstr-conv} and also the a-priori estimates \eqref{est-disc}
and \eqref{est-of-theta.omega} with $r<p$, we can see that the convergence
to 0 in \eqref{strong-hyper}, and we thus we obtain \eqref{velocity-strong}.
Using also
$\Nabla\overline u_\etau|_{I\times\varOmega}^{}\to\Nabla u_\eps|_{I\times\varOmega}^{}$
weakly* in $L^\infty(I;L^2(\varOmega;\R^d))$,
the convergence in the discrete momentum equation \eqref{ED-1+d}
is then easy.

The limit passage in the variational inequality for $\overlinemmtau$
which is behind the inclusion \eqref{ED-4+d}, cf.\ \eqref{VI}, exploits
the strong convergences  \eqref{conv-strong-+} and \eqref{conv-strong-++}.

Exploiting that $\Nabla\overlinevvtau\to\Nabla\vv_\varepsilon$ strongly in
$L_{\rm w*}^p(I;L^\infty(\varOmega;\R^{d\times d}))$, we can see that
$\overline\rr_\etau$ from \eqref{ED-4+d} converges to $\rr_\eps=\ZJ\mm_\eps$
weakly in $L^{\rexp}(I{\times}\varOmega;\R^d)$.

\medskip\noindent{\it Step 4: Convergence of the dissipation rate
and of the heat-transfer equation for $\tau{\to}0$.}
Here we use the already proved convergence in the mechanical part
together with the mechanical-energy conservation.
Using \eqref{mag-mech-engr-disc} for $l=T/\tau$,
we have the chain of estimates:
\begin{align}\nonumber
&\int_0^T\!\!\!\int_\varOmega\!\GM(\theta_\varepsilon)|\RR_\varepsilon|^2\!+
\DD|\EE(\vv_\varepsilon)|^2\!
+\NU|\Nabla\EE(\vv_\varepsilon)|^p\!+\partial_{\rr}\zeta\big(\theta_\varepsilon;
\rr_\varepsilon\big){\cdot}\rr_\varepsilon
+\varkappa|\Nabla\RR_\varepsilon|^2\,\d x\d t
\\[-.4em]\nonumber
&\quad\le\liminf_{\tau\to0}\int_0^T\!\!\!\int_\varOmega\!\!
\GM(\underline\theta_\etau)|\overlineRRtau|^2\!+
\DD|\EE(\overlinevvtau)|^2\!
\\[-.3em]\nonumber&\hspace{13em}
+\NU|\Nabla\EE(\overlinevvtau)|^p\!
+\partial_{\rr}\zeta\big(\underline\theta_\etau;\rr_\etau\big){\cdot}\rr_\etau
 +\varkappa|\Nabla\overlineRRtau|^2\,\d x\d t
 \\[-.3em]
 \nonumber&\ \,\stackrel{\eqref{mag-mech-engr-disc}}{\le}
 \!\int_\varOmega\frac\varrho2|\vv_0|^2
+\FF(\Ee_0,\mm_0)+\frac{\kappa\mu_0}{\qexp}|\Nabla\mm_0|^{\qexp}
+\frac{\mu_0}2|\Nabla u_\etau(0)|^2-\mu_0\hh_{\text{\rm geo},\tau}(0){\cdot}\mm_0
\,\d x
 \\[-.2em]\nonumber&\hspace{1em}-
 \liminf_{\tau\to0}
 \int_\varOmega\frac\varrho2|\vv_\etau(T)|^2
+\FF(\Ee_\etau(T),\mm_\etau(T))
+\frac{\kappa\mu_0}{\qexp}|\Nabla\mm_\etau(T)|^{\qexp}
+\frac{\mu_0}2|\Nabla u_\etau(T)|^2\,\d x
 \\[-.4em]\nonumber&\hspace{1em}+\lim_{\tau\to0}\bigg(
 \int_\varOmega\!\mu_0\hh_{\text{\rm geo},\tau}(T){\cdot}\mm_\etau(T)\,\d x
 +\!\int_0^T\!\!\!
 \int_\varOmega\!\varrho\overline\ff_\tau\!{\cdot}\overlinevvtau\!
 +\pdt{\hh_{\text{\rm geo},\tau}\!\!}{\cdot}\underline\mm_\etau\!
\\[-.4em]\nonumber&\hspace{5em}
-\overline\theta_\etau\CCC_\eps'(\overlinemmtau){\cdot}
\Big(\pdt{\mm_\etau\!\!}+(\overlinevvtau{\cdot}\Nabla)\overlinemmtau\Big)
-\big(\overline\theta_\etau\CCC_\eps(\overlinemmtau){+}
\phi(\overline\theta_\etau)\big){\rm div}\,\overlinevvtau\,\d x\d t\bigg)
 \\ \nonumber&\quad\le\int_\varOmega\bigg(\frac\varrho2|\vv_0|^2
+\FF(\Ee_0,\mm_0)+\frac{\kappa\mu_0}{\qexp}|\Nabla\mm_0|^{\qexp}
+\frac{\mu_0}2|\Nabla u_0|^2
-\mu_0\hh_\text{\rm geo}(0){\cdot}\mm_0-\frac\varrho2|\vv_\varepsilon(T)|^2
\\[-.4em]\nonumber&\hspace{3em}
-\FF(\Ee_\varepsilon(T),\mm_\varepsilon(T))
-\frac{\kappa\mu_0}{\qexp}|\Nabla\mm_\varepsilon(T)|^{\qexp}
-\frac{\mu_0}2|\Nabla u_\eps(T)|^2
+\mu_0\hh_{\text{\rm geo},\tau}(T){\cdot}\mm_\eps(T)\bigg)\,\d x
\\[-.2em]\nonumber&\hspace{3em}
+\int_0^T\!\!\!\int_\varOmega\!\bigg(\varrho\ff{\cdot}\vv_\varepsilon
 +\pdt{\hh_\text{\rm geo}\!\!}{\cdot}\mm_\eps\!
 -\theta_\eps\CCC_\eps'(\mm_\eps){\cdot}\DT\mm_\eps
-\big(\theta_\eps\CCC_\eps(\mm_\eps){+}
\phi(\theta_\varepsilon)\big){\rm div}\vv_\varepsilon\,\d x
\d t\bigg)
 \\[-.2em]
 &\quad= \int_0^T\!\!\!\int_\varOmega\!\!\GM(\theta_\varepsilon)|\RR_\varepsilon|^2\!+
\DD|\EE(\vv_\varepsilon)|^2\!
+\NU|\Nabla\EE(\vv_\varepsilon)|^p\!+\partial_{\rr}\zeta\big(\theta_\varepsilon;
\rr_\varepsilon\big){\cdot}\rr_\varepsilon
 +\varkappa|\Nabla\RR_\varepsilon|^2\,\d x\d t\,.
 \label{limsup-trick}\end{align}
The first inequality is due to the weak lower semicontinuity.
The last equality in \eqref{limsup-trick} is just the 
magneto-mechanical energy balance \eqref{energy+} written for the
$\varepsilon$-solution.

This magneto-mechanical energy balance follows from the tests as used
for the (formal) calculations \eqref{formula1}--\eqref{calculus-PM4}
written for the $\varepsilon$-solution. The validity of this balance 
is however not automatic and the rigorous prove needs to have granted that
testing the particular magneto-mechanical equations by $\vv_\eps$,
$\bm{S}_{\text{\sc e},\eps}$, $\RR_\eps$, $\mu_0\rr_\eps=\mu_0\ZJ\mm_\eps$, and
$\mu_0\pdt{}u_\eps$ is indeed legitimate. Here, in particular it is
important that $\bm{S}_{\rm str,\eps}{:}\Nabla\vv_\eps\in L^1(I{\times}\varOmega)$
because $\bm{S}_{\rm str,\eps}\in L^\infty(I;L^1(\varOmega;\R^{d\times d}))$
and $\Nabla\vv_\eps\in L_{\rm w*}^p(I;L^\infty(\varOmega;\R^{d\times d}))$. Also,
$$
\Delta\mm_\eps
\in\partial_{\ZJ\mm}\zeta(\theta_\eps;\ZJ\mm_\eps)
+\frac{\varphi_\mm'(\Ee_\eps,\mm_\eps){+}\CC_\mm'(\mm_\eps,\theta_\eps)}{\mu_0}
-\hh_\text{\rm geo}-\Nabla u_\eps
$$
holds pointwise a.e.\ in the sense
of $L^{\rexp'}(I{\times}\varOmega;\R^d)$, we can legitimately test it by
$\DT\mm_\eps=\rr_\eps+{\rm skw}(\Nabla\vv_\eps)\mm_\eps\in L^\rexp(I{\times}\varOmega)$ provided $r\le\min(p,\qexp^*)$ as indeed assumed, cf.\ \eqref{ass:4}.
Since $\partial_{\ZJ\mm}\zeta(\theta_\eps;\cdot)$ is single-valued except 0,
cf.\ \eqref{ass:2},
$\partial_{\ZJ\mm}\zeta(\theta_\eps;\ZJ\mm_\eps){\cdot}\ZJ\mm_\eps\in L^1(I{\times}\varOmega)$ is
single-valued.
To make the test \eqref{calculus-PM4} legitimate for the $\eps$-solution,
we need $\Nabla\pdt{}u_\eps|_{I{\times}\varOmega}\in L^1(I;L^2(\varOmega;\R^d))$. By
differentiating \eqref{Max++} in time, we have
$\Delta\pdt{}u_\eps={\rm div}(\chi_\varOmega^{}\pdt{}\mm_\eps)$. Realizing that
$\pdt{}\mm_\eps=\rr_\eps-(\vv_\eps{\cdot}\Nabla)\mm_\eps+{\rm skw}(\Nabla\vv_\eps)\mm_\eps$
belongs surely to $L^2(I{\times}\varOmega;\R^d)$, we obtain even
$\Nabla\pdt{}u_\eps|_{I{\times}\varOmega}\in L^2(I{\times}\varOmega;\R^d)$.

Altogether, this reveals that
there are actually equalities in \eqref{limsup-trick}.
Since the dissipation rate is uniformly convex in terms of rates on the
uniformly convex $L^2$- or $L^\rexp$-spaces, these rates converge not only
weakly but even strongly in these  $L^2$-spaces. Thus the dissipation rate
itself converges strongly in $L^1(I{\times}\varOmega)$.

The limit passage in the resting semilinear terms in \eqref{ED-6+d}
is then easy.

\medskip\noindent{\it Step 5: A-priori estimates uniform in $\eps>0$.}
The arguments can only slightly modify the strategy of Step~2 and the calculus
from Sect.~\ref{sec-thermodyn}. 

We exploit the magneto-mechanical energy balance as the equality used
already in \eqref{limsup-trick} yet considered on the interval $[0,t]$ and
add to it the heat-transfer equation written for the $\eps$-solution integrated
over $\varOmega{\times}[0,t]$. The nonconvexity of the stored energy 
for $\theta<\theta_{\rm c}$ now makes no problem and testing
by $\pdt{}\mm_\eps$ can be performed as in the time-continuous variant in 
\eqref{formula3}. Due to the factor $1{-}\eps$ in \eqref{ED-6+d},
we obtain the total energy balance \eqref{energy-mag} for $\eps$-solution
with ``$\eps$-part'' of the dissipative heat source ``on the left-hand
side'' like in \eqref{tot-engr-disc} which is now to be neglected for
obtaining uniform estimates with respect to $\eps>0$. Thus,
we obtain the inequality
\begin{align}\nonumber
  &\!\!\int_\varOmega\frac\varrho2|\vv_\eps(t)|^2\!
  +\varphi(\Ee_\eps(t),\mm_\eps(t))
+\frac{\!\kappa\mu_0}{\qexp}|\Nabla\mm_\eps(t)|^{\qexp}\!
+\frac{\!\mu_0}2|\Nabla u_\eps(t)|^2\!+\mu_0\hh_\text{\rm geo}(t){\cdot}\mm_\eps(t)
+\W_\eps(t)\,\d x
\\&\nonumber\hspace*{0em}
\le\int_0^t\bigg(\int_\varOmega\varrho\ff{\cdot}\vv_\eps
  +\frac{\partial \hh_\text{\rm geo}\!}{\partial t}{\cdot}\mm_\eps\,\d x
+\int_\varGamma j_{\rm ext\,}\d S\bigg)\,\d t
\\[-.1em]&\hspace*{2.5em}
+\int_\varOmega\frac\varrho2|\vv_0|^2\!+\varphi(\Ee_0,\mm_0)
+\frac{\kappa\mu_0}{\qexp}|\Nabla\mm_0|^{\qexp}\!+\frac{\mu_0}2|\Nabla u_0|^2\!
+\mu_0\hh_\text{\rm geo}(0){\cdot}\mm_0+\W_0\,\d x\,.
\label{energy-mag-eps}\end{align}

By using the Young and the Gronwall inequalities, we obtain the a-priori
estimates uniform with respect to $\eps>0$:
\begin{subequations}\label{est-esp}\begin{align}\label{est-esp1}
    &\|\vv_\eps^{}\|_{L^\infty(I;L^2(\varOmega;\R^d))}^{}\le C\,,
  \\\label{est-esp5}
&\|\Ee_\eps^{}\|_{L^\infty(I;L^2(\varOmega;\R_{\rm sym}^{d\times d}))}^{}\le C\,,
   \\\label{est-esp4}
   &\|\mm_\eps^{}\|_{L^\infty(I; H^1(\varOmega;\R^d)) }^{}\le C\,,
   \\\label{est-esp7}
   &\|u_\eps^{}\|_{L^\infty(I;H^1(\R^d))}^{}\le C\,,
 \\\label{est-esp-w}
   &\|\W_\eps^{}\|_{L^\infty(I;L^1(\varOmega))}^{}\le C\ \ \text{ and }\ \
   \|\theta_\eps^{}\|_{L^\infty(I;L^1(\varOmega))}^{}\le C\,.
 \end{align}\end{subequations}

Then, we employ the magneto-mechanical energy balance itself as the equality 
used already in \eqref{limsup-trick} for the $\eps$-solution.
Now, in the time-continuous situation, we can use
also the identity $\rr_\varepsilon=\ZJ\mm_\varepsilon$.
Thus, as in \eqref{mag-mech-engr-disc+}, we arrive to:
\begin{align}\nonumber
&
\int_0^t\!\!\int_\varOmega\!\!\GM(\theta_\varepsilon)|\ZJ\Ep_\varepsilon|^2\!+
\DD|\EE(\vv_\varepsilon)|^2\!
+\NU|\Nabla\EE(\vv_\varepsilon)|^p\!+\partial_{\ZJ\mm}\zeta\big(\theta_\varepsilon;
\ZJ\mm_\varepsilon\big){\cdot}\ZJ\mm_\varepsilon
\\[-.5em]\nonumber&\hspace{22em}
+\varkappa|\Nabla\ZJ\Ep_\varepsilon|^2+\delta_s|\theta_\etau|^s
 +\delta_\sigma|\Nabla\theta_\etau|^\sigma\,\d x\d t
\\[-.2em]\nonumber&\hspace{.1em}
=\int_\varOmega
\frac\varrho2|\vv_0|^2\!
+\varphi(\Ee_0,\mm_0)+\frac{\kappa\mu_0}{\qexp}|\Nabla\mm_0|^{\qexp}\!
+\frac{\mu_0}2|\Nabla u_0|^2\!-\mu_0\hh_\text{\rm geo}(0){\cdot}\mm_0\!
+\mu_0\hh_\text{\rm geo}(T){\cdot}\mm_\eps(T)\,\d x
\\[-.1em]
&\hspace{1em}
+\int_0^t\!\!
\int_\varOmega\!2\varrho\ff{\cdot}\vv_\varepsilon
 +2\pdt{\hh_\text{\rm geo}\!\!}{\cdot}\mm_\eps\!
+\big|\theta_\eps\CCC_\eps'(\mm_\eps){\cdot}\DT\mm_\eps\big|
+\big|\big(\theta_\eps\CCC_\eps(\mm_\eps){+}
\phi(\theta_\varepsilon)\big){\rm div}\vv_\varepsilon\big|\,\d x
\d t \,.
 \label{limsup-trick+}\end{align}

Exploiting \eqref{est-esp-w} with the estimation as
\eqref{adiabat-est-1}--\eqref{adiabat-est-2}, we can treat the
adiabatic terms on the right-hand side of \eqref{limsup-trick+}.
We thus obtain the estimates
\begin{subequations}\label{est-eps}\begin{align}\label{est-eps1}
    &\|\vv_\eps^{}\|_{ L^p(I;W^{2,p}(\varOmega;\R^d))}^{}\le C\,,
\\&\|\ZJ\Ep_\varepsilon\|_{L^2(I;H^1(\varOmega;\R^{d\times d}))}^{}\le C\,,
\\&\|\ZJ\rr_\varepsilon\|_{L^{\rexp}(I{\times}\varOmega;\R^d))}^{}\le C\,,
\\&\|\theta_\eps\|_{L^\sigma(I;W^{1,\sigma}(\varOmega))
  \,\cap\,L^s(I{\times}\varOmega)}^{}\le C_{s,\sigma}\,,
\\[-.4em]\label{est-disc+vartheta+}
&\|\W_\eps\|_{L^\sigma(I;W^{1,\sigma}(\varOmega))
 \,\cap\,L^s(I{\times}\varOmega)}^{}\le C_{s,\sigma}
\ \ \text{ with $\ \ 1\le s<1+\frac2d$,\
$\ 1\le\sigma<\frac{d{+}2}{d{+}1}$.}
\end{align}\end{subequations}
The estimates \eqref{est-disc9}, \eqref{est-Delta}, and \eqref{Gronwall-3-2}
for $\eps$-solution uniform with respect to $\eps>0$ can be obtained
analogously as above.

\medskip\noindent{\it Step 6: Convergence for $\eps\to0$.}
This final convergence towards a weak solution due to Definition~\ref{def}
copies the arguments in the Step~4 above. By the Banach selection principle,
we obtain a subsequence converging weakly* with respect to the topologies
indicated in estimates in Step~5
to some limit $(\vv,\Ee,\Ep,\RR,\mm,u,\theta)$. The only small difference is
that (\ref{conv-strong}f,g) and (\ref{conv-strong+}a,b)  now involves also
the convergence $\CCC_\eps\to\CCC$,
$\CCC_\eps'\to\CCC'$,
$\CC_\eps\to\CC$, and $[\CC_\eps]_\mm'\to\CC_\mm'$.

{
\section*{{\large Acknowledgments}}

\vspace*{-1em}

This research has been partially supported also from the CSF (Czech Science
Foundation) project 19-04956S, from the M\v SMT \v CR (Ministry of Education of
the Czech Rep.) project CZ.02.1.01/ 0.0/0.0/15-003/0000493, and by the
institutional support RVO: 61388998 (\v CR).

\medskip
}

\end{document}